\documentclass[9pt,shortpaper,twoside,web]{ieeecolor}
\usepackage{generic}
\usepackage{cite}
\usepackage{amsmath,amssymb,amsfonts}
\usepackage{algorithmic}
\usepackage{graphicx}
\usepackage{textcomp}

\usepackage{amsopn, amstext}
\usepackage{colortbl}
\usepackage{empheq}
\usepackage{graphicx}
\usepackage{color}
\usepackage[OT1]{fontenc}
\usepackage{mathtools}
\usepackage{scalerel}
\usepackage{float}

\usepackage{enumerate}

\newtheorem{remark}{Remark}
\newtheorem{theorem}{Theorem}

\newtheorem{lemma}{Lemma}

\newtheorem{corollary}{Corollary}
\newtheorem{proposition}{Proposition}

\def\s{\small}
\def\n{\normalsize}

\def\+{\!+\!}
\def\-{\!-\!}
\def\={\!=\!}

\def\E{\mathbb{E}}
\def\Eo{\E^{W_0}}
\def\x{\bar{x}}

\def\u{\bar{u}}

\def\mD{\mathcal{D}}

\def\E{\mathbb{E}}

\def\S{\mathbb{S}}

\def\I{\mathcal{I}}

\def\A{\mathcal{A}}

\def\B{\mathcal{B}}

\def\C{\mathcal{C}}

\def\U{\mathcal{U}}

\def\J{\mathcal{J}}

\def\F{\mathcal{F}}

\def\H{\mathcal{H}}

\def\bF{\mathbb{F}}

\def\R{\mathbb{R}}

\def\xN{x^{(N)}}

\def\uN{u^{(N)}}

\def\bx{\mathbf{x}}

\def\bu{\mathbf{u}}

\def\+{\!+\!}
\def\-{\!-\!}
\def\={\!=\!}
\def\x{\bar{x}}

\def\u{\bar{u}}

\def\s{\small}

\def\n{\normalsize}

\def\T{\int_{0}^{T}}

\def\Eo{\mathbb{E}_0}

\def\hmc{{\text{\rm\textbf{H\tiny MC}}}}

\def\hmg{{\text{\rm\textbf{H\tiny MG}}}}

\def\hmt{{\text{\rm\textbf{H\tiny MT}}}}

\def\ccmg{{\text{\rm\textbf{CC\tiny MG}}}}

\def\ccmt{{\text{\rm\textbf{CC\tiny MT}}}}

\def\1mm{\vspace{-1mm}}
\def\2mm{\vspace{-2mm}}
\def\3mm{\vspace{-3mm}}
\def\4mm{\vspace{-4mm}}

\def\BibTeX{{\rm B\kern-.05em{\sc i\kern-.025em b}\kern-.08em
    T\kern-.1667em\lower.7ex\hbox{E}\kern-.125emX}}
\begin{document}
\title{A Unified Relation Analysis of Linear-quadratic Mean-field Game, Team and Control}
\author{Jianhui Huang, Zhenghong Qiu, Shujun Wang, and Zhen Wu
\thanks{This work was supported in part by RGC 153275/16P, P0031141, P0030808, P0008686, P0031044, the National Natural Science Foundation of China (No. 11831010, 61961160732, 12001320), the Natural Science Foundation of Shandong Province (No. ZR2019ZD42), the Taishan Scholars Climbing Program of Shandong (No. TSPD20210302) and the Young Scholars Program of Shandong University. The authors also acknowledge the financial support from PolyU-SDU joint research center.  \emph{(Corresponding author: Shujun Wang.)}}
\thanks{Jianhui Huang is with the Department of Applied Mathematics, The Hong Kong Polytechnic University, Hong Kong (e-mail: majhuang@polyu.edu.hk). }
\thanks{Zhenghong Qiu is with the Department of Applied Mathematics, The Hong Kong Polytechnic University, Hong Kong (e-mail: zhenghong.qiu@connect.polyu.hk). }
\thanks{Shujun Wang is with the School of Management, Shandong University, Jinan, Shandong 250100, China (e-mail: wangshujun@sdu.edu.cn).}
\thanks{Zhen Wu is with the School of Mathematics, Shandong University, Jinan, Shandong 250100, China (e-mail: wuzhen@sdu.edu.cn).}}

\maketitle

\begin{abstract}
This paper revisits well-studied dynamic decisions of weakly coupled large-population (LP) systems.
Specifically, three types of LP decision problems: mean-field game (MG), mean-field team (MT), and mean-field-type control (MC), are completely analyzed in a general stochastic linear-quadratic setting with controlled-diffusion in state dynamics and indefinite weight in cost functional. More importantly, interrelations among MG, MT and MC are systematically discussed; some relevant interesting findings are reported that may be applied to a structural analysis of general LP decisions.
\end{abstract}


\section{Introduction} \label{sec1}
Purpose of this paper is to investigate possible connections among mean-field type control (MC) problem, mean-field game (MG) problem, and mean-field social optima problem (also named as mean-field team (MT)). All three problems have been well investigated along their individual research lines from various perspectives in the last few decades. One common feature shared or implied by the modelings of these three problems, is the existence of considerable interactive agents or decision makers in underlying dynamic systems. Usually, such systems are also recognized as large-population (LP) system (or multi-agent system, in a more general sense) to emphasize complex interactions among a large number of coupled agents. More details of LP systems are presented as follows.

\subsection{LP system}

LP systems have been widely applied in a wide variety of areas, including economics, biology, engineering, and social science (see, \cite{bensoussan2013mean,HCM2015,weintraub2008markov}, etc). The most salient
feature of LP systems is the existence of sufficiently many agents whose dynamics and cost functionals are interrelated via some empirical state-average or control-average. Although the effect of an individual agent is negligible in a micro-scale, but whole effects of their
statistical behaviors, cannot be ignored in a macro-scale. Last few decades have witnessed a rapid growth in research for dynamic decisions of LP systems, possibly because of burgeoning practical models with large-scaled interactions arising in the real world.

When the number of agents is sufficiently large, there will emerge two severe issues in analysis of dynamic decisions of LP systems: (i) The system complexity, due to structural couplings, will dramatically increase along with the growth of the agent's population size; so, it becomes unrealistic for a given agent in these systems to compile all other agents' information (centralizing information) to tackle related highly coupled problems. (ii) Due to the coupling structure, possible analysis on LP systems is intrinsically of high-dimensionality and becomes vulnerable to the so-called ``curse of dimensionality'' (e.g., \cite{r1997,ks2005}). Indeed, computational complexity would increase exponentially in practical numerical analysis.

Subsequently, mean-field method has drawn increasing research attention because it
provides an effective scheme to obtain an asymptotically optimal
decentralized control based on decentralized information of individual
agent with much lower computational burden. In particular, mean-field
method enables each agent to take advantage of mean field interaction to
convert analysis of the original LP problem (centralizing and
high-dimensional) to an optimization problem (decentralized
and low-dimensional) parameterized by a sufficient statistics of aggregation effects of other
agents. Interested readers are referred to
\cite{cardaliaguet2018mean, slanina2001mean}
for mean-field method in  economics, \cite{ma2011decentralized, couillet2012electrical} in
engineering, \cite{bauch2004vaccination, breban2007mean} in biology,
as well as \cite{graber2016linear, huang2007large} in management science or
operational research.

\subsection{MG problem}
The above mean-field method has been extensively applied in game models (see  \cite{bensoussan2013mean, HCM2015, carmona2018prob, huang2007large, lasry2007mean}). Agents in these aforementioned
works are non-cooperative because they aim to optimize competitive objectives, thus relevant models are named as mean-field game (MG). Introduction of MG can be traced back to parallel pioneering works of Lasry and
Lions \cite{lasry2007mean} from a partial differential equation (PDE) perspective, and of Huang, Caines, and Malham\'e \cite{huang2007large} from an engineering perspective.  Specifically, in \cite{huang2007large} an $\varepsilon$-Nash equilibrium is designed via Nash certainty equivalence (NCE) approach which is
also called the consistency condition (CC).

In \cite{lasry2007mean}, MG is studied via two limiting coupled partial differential equation systems which are inspired  from physical
particle systems. The first one is Hamilton-Jacobi-Bellman equation which describes the reaction of agents to population aggregations. The second one is Fokker-Planck-Kolmogorov equation which describes behaviors of the population aggregation. For a more detailed discussion of MG, interested readers may refer \cite{bensoussan2016linear},
  \cite{hu2018linear}, \cite{HL2018}, \cite{huang2016backward} for linear-quadratic-Gaussian (LQG)
  MG; \cite{carmona2013probabilistic} for probabilistic analysis of MG; \cite{moon2019risk} for risk-sensitive type MG;
  \cite{moon2014discrete} for some discrete-time MG.

\subsection{MT problem}
Apart from noncooperative MG,
cooperative multi-agent decision problems (social optima) are also research hot topics due to their
theoretical interests and real application potentials. In MT problems, all the
agents usually work cooperatively to minimize a social cost as the sum of $N$ individual costs
containing mean-field coupling.

The recent research of MT  largely adopt two different methods. The first method may be termed as fixed-point approach that includes mainly three steps. In first step, each individual agent starts by constructing an auxiliary control problem by applying variational method, duality and mean-field approximations to the social cost functional. Then, second step aims to solve the auxiliary control problem for which the relevant auxiliary optimal control involving some pre-frozen mean-field terms. So, the last third step would determine such mean-field terms via some consistency condition analysis. Herein, some fixed point analysis will be involved to achieve such consistency hence the approach is termed as fixed-point one. More details of the fixed-point approach can be referred to \cite{carmona2018prob,lz2008,huang2007large,HCM2012}, etc.

Another method is called direct approach that is free of the fixed-point argument. The rough idea of direct approach is as follows. As the first step, each agent starts by formally solving the MT problem as a high dimensional control problem in a direct manner. Then, an optimal control can be represented in a feedback form via Riccati-type equations that should converge to some standard but low-dimensional Riccati equations. The, the next step is derive the limit of the optimal control through limiting Riccati equations when the population size tends to infinity. Here, some block structure analysis will be applied to reduce the dimension. For more discussion of direct approach method, readers are referred to
\cite{huangmeng2019linear, wang2020indefinite}, etc.

Till now, there have accumulated very rick literature on MT studies. For a more comprehensive review, interested readers may
refer \cite{HCM2012} for centralized and decentralized controls of a LQG MT problem where the
asymptotic optimality of mean-field controls is also illustrated;
\cite{nourian2012nash} for the mean-field social solution of consensus problems;
\cite{wang2017social} for social optima of mean-field LQG control models with Markov jump parameters; \cite{huang2019linear} for an LQG MT problem involving a major player.

\subsection{MC problem}
In the aforementioned MT problem, each agent has the freedom to choose its own (decentralized) strategy but are subject to different information sets that may vary from agent to agent. On the other hand, for a system consisting of interacting agents who are framed with an identical and centralized information set, some MC type control problems will be suggested when we consider relevant dynamic decisions (see \cite{andersson2011maximum, bardi2014linear, BLP2014, NNY2020, priuli2015linear}, etc). Such MC problem shares
a similar form as proposed in MT problem, and the mean-field term is now uniformly
determined by a centralizing system instead of being affected by the aggregation of the population. MC problems aim to design a strategy to all agents at the same time, such that the resulting behavior is optimal with respect to costs imposed by a centralized planner. In MC
problems, the mean-field
term is influenced by the agent individually. Indeed, the
state equation also contains the probability distribution of the state and thus  is called the
McKean-Vlasov stochastic differential equation (SDE), a type of mean-field forward stochastic differential equation
(MF-FSDE), which was suggested by Kac \cite{kac1956foundations} in 1956 as a stochastic model for
the Vlasov
kinetic equation of plasma, and a study of this type of equation was initiated by McKean \cite{mckean1966class} in
1966.
\subsection{The analysis outline and comparison}

Noting that the three problems MG, MT, MC share a similar mathematical form, so a natural question is:
is there any intrinsic relations among them? In \cite{graber2016linear}, the relation between MG and MC has
been studied. The author firstly derived the feedback form optimal control of the MC problem by
using  decoupling method consisting of two symmetric Riccati equations. Further,  it was proved that
 the MC optimal control is also the equilibrium strategy for a class of MG problems in certain
 cases. Last the asymptotic optimality of MG strategies in $N$-player game has also been studied.
 Motivated by \cite{graber2016linear}, this paper aims to further study the relations among MG, MT
 and MC problems. Specifically, the MG, MT and MC problems are studied under a unified mathematical
 framework. Specifically, empirical average terms of both state and control enter the cost functionals and state dynamics
 (in both drift and diffusion terms).

To conclude, main contributions of this paper can be summarized as follows:
\begin{enumerate}
  \item We investigate the three problems of MG, MT and MC, in a unified framework where the state dynamics takes a quite general form. Specifically, its diffusion term depends on control and control-average both, and is driven by individual noises and a common noise at the same time. In control literature, related stochastic
LQG problem can be referred to have multiplicative-noise. Moreover, the inclusion of control variable in diffusion is well motivated by various real applications. One such example is the well-known mean-variance portfolio problem \cite{37,41} where the control (risky portfolio allocation) process naturally enters the diffusion of wealth process. Our general setting also makes our research to be an extension of many existing literature \cite{graber2016linear, HCM2012, wang2020indefinite, yong2013linear}, etc.
   \item We study relations among the optimal control of MC and the mean-field controls of MG and MT problem. Through analyzing the Hamiltonian system of MC problem  and  the CC systems of MG and MT problem, we find that the optimal control of MC problem is equivalent to the mean-field control of MT problem. Such result provides short-cut to deal with MT problem in practical application. Each agent in the system could directly calculate its own related MC problem instead of deriving an auxiliary control problem (fixed point approach) or computing the limit of the centralized optimal control (direct approach). Moreover, we also derive that in certain cases, the optimal control of MC problem is also the mean-field control of MG problem. Such result is consist to that in \cite{graber2016linear} which can be treated as a special case of the current study.
       %Last, we compare the mean-field controls of MT problem derived by fixed point approach and direct approach and also find that they are equivalent. These two routes leads to a same mean-field control.
  \item We study relations among value functions of MG, MT and MC. By some asymptotic analysis and SDE estimations, it is shown that in general value functions of MT and MC is less than that of MG, while value functions of MT and MC are very close to each other (specifically, their difference is of $O\left(N^{-\frac{1}{2}}\right)$ order). Moreover, in some special cases value functions of MG, MT and MC could be all close together (in $O\left(N^{-\frac{1}{2}}\right)$ order).
\end{enumerate}

%The rest of this paper is organized as follows. Section \ref{sec 2} introduces the basic structure of the LP system. The relations among the optimal controls of MC and optimal strategies of MG and MT problem are analyzed, and some numerical analysis is given to illustrate the theoretical results in Section \ref{sec 4}. Section \ref{sec 5} concludes this work.

This paper is organized as follows. Section \ref{sec 2} introduces the basic structure of the LP system. In Section \ref{sec 4}, relations among the optimal control of MC and the optimal strategies of MG and MT problem are analyzed, and some numerical analysis is given to illustrate theoretical results. Section \ref{sec 5} concludes this work.

\quad\\
\emph{\textbf{Notation}}: Throughout this paper, %$\prod$ denotes the Cartesian product. For a set $A$, $|A|$ denotes its cardinality.
$\mathbb{R}^{n\times m}$ and
$\mathbb{S}^n$ denote the sets of all $(n\times m)$ real matrices and
 all $(n\times n)$ symmetric matrices respectively. If $v$ consists of several sub-vectors, $v_1,\cdots,v_N$, it is denoted for simplicity of notation as $(v_1,\cdots, v_N)$ instead of $((v_1)^{\top},\cdots,(v_N)^{\top})^{\top}$, where $v^{\top}$ is the transpose of $v$. $I$ is
the identity matrix. We denote $\|\cdot\|$ as the standard Euclidean
norm and $\langle\cdot,\cdot\rangle$ as the standard Euclidean inner
product. For a vector $v$ and a symmetric matrix $S$, $\|v\|_S^2:=\langle Sv,v\rangle=v^{\top}Sv$. %For a linear subspace $\Lambda$, $v\bot\Lambda$ denotes that vector $v$ is vertical to $\Lambda$ (i.e., for any vector $v'\in\Lambda$, $\left\langle v, v'\right\rangle=0$).
%For
%a matrix $M$, the norm $\|M\|=\sqrt{\text{Tr}(M^{\top}M)}$, and the
%max-norm, which is equivalent to the maximum absolute value of all
%elements, is denoted by $\|M\|_{\max}$.
 $S> 0 \ (\geq 0)$ means that
$S$ is positive (semi-positive) definite, and $S\gg 0$ means that, $S
- \varepsilon I \geq 0$, for some $\varepsilon>0$. %For a matrix $M$, $\lambda_{\max}(M)$ denotes the maximum eigenvalue of $M$, while $\lambda_{\min}(M)$ denotes the minimum eigenvalue of matrix $M$.

%We denote $\mathcal F^{i}_t$ as the filtration generated by the $i^{\text{th}}$ Brownian motion $W_i(t)$: $\sigma\{W_i(s), 0\leq s\leq t\}$. For a Euclidean space $\mathbb{H}$, let $\varsigma\in L^2_{\mathcal F_T}(\Omega;\mathbb{H})$ be the set of $\mathbb{H}$-valued $\mathcal F_T$ measurable random variable such that $\E |\varsigma|^{2}< + \infty$ and $f\in L^\infty_{\F}(s,T;\mathbb{H})$ be the set of $\{\mathcal F_t\}_{t\geq 0}$-progressively measurable, $[s,T]\!\times\!\Omega\rightarrow\mathbb{H}$ stochastic process and the essential supremum $\text{esssup}_{(t,\omega)\in[s,T]\times\Omega}|f(t)|<\infty$.
%$f\in L^p(s,T;\mathbb{H})$ denotes the set of $\mathbb{H}$-valued function such that $\int_{s}^{MT}|f(t)|^pdt<\infty$ and $f\in L^p_{\F}(s,T;\mathbb{H})$ denotes the set of $\{\mathcal F_t\}_{t\geq 0}$-progressively measurable, $[s,T]\times\Omega\rightarrow\mathbb{H}$ stochastic process such that $\E  \int^{\top}_0|f(t)|^pdt<\infty$.

%Throughout this paper,
We suppose that $(\Omega,\F,\mathbb{P})$ is an complete probability space, and $W(\cdot)$ = $(W_{0}(\cdot)$,$W_{1}(\cdot)$,$\cdots$,$W_N(\cdot))$ is a $(N+1)$-dimensional standard Brownian motion defined on it. Here, $W_0(\cdot)$ stands for the common noise which are assumed to be seen by all agents; while $W_i(\cdot)$ stands for the individual information observed by $i^{th}$ agent, but cannot be observed by other agents.
Denote $\sigma$-algebra $\F(t):=\sigma\{W_i(s), 0\leq s\leq t,0\leq i\leq N\}$ and $\bF:=\{\F(t)\}_{t\geq0}$ corresponding to $W(\cdot)$.
$\Eo$ denotes the conditional expectation w.r.t. the filtration generated by $W_0(\cdot)$.
%
%% $\mathcal{G}_t^i=\F_t^i\bigvee\sigma\{\xi_i\}$ for $1\leq i\leq N$, $\mathcal{G}_t=\F_t\bigvee\sigma\{\xi_i,1\leq i\leq N\}$, $\mathcal{H}_t^i=\sigma\{x_s^i,0\leq s\leq t\}$ for $1\leq i\leq N$, $\mathcal{H}_t=\sigma\{x_s^i,0\leq s\leq t,1\leq i\leq N\}$.
%  Correspondingly, we denote filtration
%% $\mathbb{G}^i=\{\mathcal{G}_t^i\}_{0\leq t\leq T}$ and $\mathbb{H}^i=\{\mathcal{H}_t^i\}_{0\leq t\leq T}$ for $1\leq i\leq N$,
% $\mathbb{F}=\{\F_t\}_{0\leq t\leq T}$.
%% $\mathbb{G}=\{\mathcal{G}_t\}_{0\leq t\leq T}$.
 Next, for a given Euclidean space $\mathbb{H}$ and filtration $\mathbb{G}$, we introduce the following spaces:\s
\begin{equation}
\begin{aligned}
%&L^{2}(0,T;\mathbb{H})=\big\{x : [0,T]\rightarrow \mathbb{H}\big| \    {\T}\|x(t)\|^2dt<\infty\big\},\\
&L^{\infty}(0,T;\mathbb{H})=\big\{x : [0,T]\rightarrow \mathbb{H}\big| \ x(\cdot) \text{ is bounded and measurable}\big\},\\
& C([0,T];\mathbb{H})=\big\{x : [0,T]\rightarrow \mathbb{H}\big| \ x(\cdot) \text{ is continuous}\big\},\\
&L_{\mathbb{G}}^2(0,T;\mathbb{H})=\big\{x : [0,T]\times\Omega\rightarrow \mathbb{H}\big| \ x(\cdot) \text{ is }\mathbb{G}\text{-progressively}\\
&\hspace{2.3cm}\text{measurable, }\|x(t)\|_{L^2}^2:=\E  {\T}\|x(t)\|^2dt<\infty\big\}.
\end{aligned}
\end{equation}

\section{Problem formulation}\label{sec 2}

In this paper, on $(\Omega,\F,\bF,\mathbb{P})$ we consider an LP system with $N$ agents, denoted by $\{\A_{i}\}_{i\in\I}$, and $\I=\{1,\cdots,N\}$ denotes the index set of the agents. The aggregation of all  agents is denoted by $\A:=\{\A_{i}\}_{i\in\I}$. The state process of the $i^{\text{th}}$ agent $ \A_i$ is modeled by a controlled linear SDE on finite time horizon $[0,T]$\s
\begin{equation}\label{xN}
\mathcal{S}_i\left\{\begin{aligned}
& dx_i \= (Ax_i \+ \bar{A}\xN  \+ Bu_i \+ \bar{B}\uN )dt\\
&\hspace{1cm} \+ (Cx_i \+ \bar{C}\xN  \+ Du_i \+ \bar{D}\uN)dW_i(t)\\
& \hspace{1cm}  \+  (C_0x_i \+ \bar{C}_0\xN \+ D_0u_i \+ \bar{D}_0\uN)dW_0(t),\\
& x_i(0) = x_0,
\end{aligned}\right.
\end{equation}\n
where $\xN = \frac{1}{N}{\sum_{i\in\I}x_i}$, $\uN = \frac{1}{N} {\sum_{i\in\I}u_i} $, and for the sake of notation simplicity, we denote  $\bu := (u_1,\cdots,u_N)$ and $\bx := (x_1,\cdots,x_N)$. Then $\mathcal{S} := \{\mathcal{S}_i\}_{i\in\I}$ forms a weakly coupled large population  system, since each agent $\A_i$ is coupled with the others via $\xN$ and $\uN$, and  $\A_i$ could only provide weak influence to the others  at $O\left(N^{-1}\right)$ order. In system $\mathcal{S}$, by
such interactive coupling, we can introduce the following  information structures
\begin{description}
  \item[(i)]   \emph{Centralized, open-loop information}: The centralized, open-loop information is $\bF:=\{\F(t)\}_{t\geq0}$, which is generated by all the Brownian motions $\{W_i\}_{i=0}^{N}$. It is the largest information structure.
  \item[(ii)]  \emph{Centralized, closed-loop information}: Let $\H_i(t) = \sigma\{x_i(s),0\leq s\leq t\}$, $\H(t) = \H_1(t)\vee\cdots\vee\H_N(t)$ and $\mathbb{H} = \{\H(t)\}_{t\geq0}$. Then the centralized, closed-loop information is $\mathbb{H}$, which is generated by  the information  of all the state processes $\{x_i\}_{i\in\I}$.
  \item[(iii)]  \emph{Decentralized, open-loop information}:  Let $\F_i(t) = \sigma\{W_0(s),W_i(s),0\leq s\leq t\}$  and $\bF_i = \{\F_i(t)\}_{t\geq0}$. Then the decentralized, open-loop information is $\bF_i$, which is generated by  the information  of the $i^{\text{th}}$ Brownian motions $W_i$.
  \item[(iv)] \emph{Decentralized, closed-loop information}: Let $\mathbb{H}_i = \{\H_i(t)\}_{t\geq0}$. Then the decentralized, closed-loop information is $\mathbb{H}_i $, which is generated by  the information  of the $i^{\text{th}}$ state process $x_i$.
\end{description}
\begin{remark}
The relations of the information structures are given as
  \begin{equation*}
   \begin{aligned}
  &\H_i(t)\subset\H(t)\subset\F(t),\quad \F_i(t)\subset\F(t).
  \end{aligned}
  \end{equation*}
\end{remark}
Based on such aforementioned information structures, we can introduce the following admissible control sets
\begin{description}
  \item[(i)]   \emph{Centralized, open-loop admissible control set}:
  \begin{equation*}
  \begin{aligned}
  \U_c^{ol} = \Big\{u|u\in L_{\mathbb{F}}^2(0,T;\R^m)\Big\}.
  \end{aligned}
  \end{equation*}
  \item[(ii)]  \emph{Centralized, closed-loop admissible control set}:
  \begin{equation*}
  \begin{aligned}
  \U_c^{cl} = \Big\{u|u\in L_{\mathbb{H}}^2(0,T;\R^m)\Big\}.
  \end{aligned}
  \end{equation*}
  \item[(iii)]  \emph{Decentralized, open-loop admissible control set}:
  \begin{equation*}
  \begin{aligned}
  \U_i^{ol} = \Big\{u|u\in L_{\mathbb{F}_i}^2(0,T;\R^m)\Big\},\quad i \in \I.
  \end{aligned}
  \end{equation*}
  \item[(iv)] \emph{Decentralized, closed-loop admissible control set}:
  \begin{equation*}
  \begin{aligned}
  \U_i^{cl} = \Big\{u|u\in L_{\mathbb{H}_i}^2(0,T;\R^m)\Big\},\quad i \in \I.
  \end{aligned}
  \end{equation*}
\end{description}
\begin{remark}
  The admissible control sets involve the following  relations
  \begin{equation*}
   \begin{aligned}
  &\U_i^{cl}\subset\U_c^{cl}\subset\U_c^{ol},\quad \U_i^{ol}\subset\U_c^{ol}.
  \end{aligned}
  \end{equation*}
\end{remark}
To evaluate each control law $u_i$, we introduce the individual cost functional as\s
\begin{equation}\label{JMG}
\begin{aligned}
&\J_i(x_0;\bu)\=\J_i(x_0;{u}_i,{\bu}_{-i})\=\frac{1}{2}\E \Bigg\{{\T}\Big[\|x_i\-\Gamma_1\xN
\|_{Q}^{2}\\
&\hspace{2.5cm}\+ \|u_i\|^{2}_{R}\Big]dt +
\|x_i(T)\-\Gamma_2\xN(T)\|_{G}^{2}\Bigg\},\\
\end{aligned}
\end{equation}\n
where  $\bu_{-i} = (u_1,\cdots,u_{i-1},u_{i+1},\cdots,  u_N)$. By letting $\J_i^{\textrm{\tiny MG}} = \J_i$, we propose the following MG problem

\quad\\
(\textbf{MG}): For given initial value $x_0\in\R^n$, each agent $\A_i$ finds a decentralized open-loop strategy $\u_i\in \U_i^{ol}$ such that
\begin{equation}\label{MG 2}
\begin{aligned}
\J_i^{\textrm{\tiny MG}}(x_0;\u_i,\bar{\bu}_{-i}) = \inf_{u_i \in \U_i^{ol}}\J_i^{\textrm{\tiny MG}}(x_0;u_i,\bar{\bu}_{-i}),
\end{aligned}
\end{equation}
 where $(\bx,\bu)$ is subject to $\mathcal{S}$.

\quad\\
Any $\bar{\bu} = (\u_1,\cdots,\u_N)$ satisfying \eqref{MG 2} is called an optimal strategy set and the
corresponding state process set $\bar{\bx} = (\x_1(x_0;\bar{\bu}),\cdots,\x_N(x_0;\bar{\bu}))$ is called an optimal state process set. The pair $(\bar{\bx}, \bar{\bu})$ is called an optimal pair.

The agents in the aforementioned problem (\textbf{MG}) are competitive. Each agent aims to minimize its own cost functional against the others. By contrast, when we study a cooperative system, it suffices to consider the social cost functional as
\begin{equation}\label{Jsoc}
\begin{aligned}
\J^{\textrm{\tiny \rm MT}}(x_0;{\bu}) = \sum_{i\in\I}\J_i(x_0;{\bu}).
\end{aligned}
\end{equation}
Then we pose the MT problem.

\quad\\
(\textbf{MT}): For given initial value $x_0\in\R^n$, all the agents $\A_1,\cdots,\A_N$ cooperate to find a decentralized open-loop strategy set $\bar{\bu}=(\u_1,\cdots,\u_N )\in \U_1^{ol}\times\cdots\times \U_N^{ol}$   such that \[\J^{\textrm{\tiny \rm MT}}(x_0;\bar{\bu}) = \inf_{\bu \in \U_1^{ol}\times\cdots\times \U_N^{ol}}\J^{\textrm{\tiny \rm MT}}(x_0;\bu),\]
 where $(\bx,\bu)$ is subject to $\mathcal{S}$.

\quad\\
Similarly we have the optimal strategy set $\bar{\bu}$, optimal state process set $\bar{\bx}$ and the optimal pair $(\bar{\bx}, \bar{\bu})$.
\begin{remark}
  Observe that in (\textbf{MG}) and (\textbf{MT}), each agent chooses its control in decentralized
  admissible control set   $\U_i^{ol}$.  This is because if we consider the centralized admissible control set   $\U_c^{ol}$, it will face some
  difficulties in the practical application. Firstly, the agent may be only
  able to access its own information, and the information of the others may be unavailable for it in some practical models (see \cite{aoki1972feedback,lau1972decentralized,wang1973stabilization}). Secondly, by the coupling structure, the dynamic optimization will be
subjected to the curse of dimensionality and complexity in numerical analysis
 in practice (see \cite{r1997,ks2005}). Thus,  to some extent,  decentralized
  control would be more practicable in real application than centralized
  control (see \cite{huang2007large,km2014}).
\end{remark}

Note that in \eqref{xN} and \eqref{Jsoc}, all the agents share the same object: to minimize $\J^{\textrm{\tiny \rm MT}}$, and their dynamics \eqref{xN} are exchangeable. Thus, it is sufficient to make a heuristic reasoning that all the agents would act  identically to achieve the social optima. Moreover, by observing the weak coupling structure in dynamics \eqref{xN}, each agent could only provide insignificant influence to the others  at $O\left(N^{-1}\right)$ order.  Consequently, we could also conjecture that all the agents are mutually independent when $N\rightarrow \infty$, and in this case we have $\lim_{N\to \infty}\xN=\Eo x_1 = \cdots = \Eo x_N$ and $\lim_{N\to \infty}\uN=\Eo u_1 = \cdots = \Eo u_N$ by the exchangeability. Thus, it suffices to consider the following limiting dynamic
\begin{equation}\label{E x}
\left\{\begin{aligned}
& dx_i \= (Ax_i \+ \bar{A}\Eo x_i \+ Bu_i \+ \bar{B}\Eo u_i )dt\\
& \hspace{1cm} \+ (Cx_i \+ \bar{C}\Eo x_i \+ Du_i \+ \bar{D}\Eo u_i)dW_i\\
& \hspace{1cm}  \+  (C_0x_i \+ \bar{C}_0\Eo x_i \+ D_0u_i \+ \bar{D}_0\Eo u_i)dW_0,\\
& x_i(0) = x_0,
\end{aligned}\right.
\end{equation}
and limiting cost functional
\begin{equation}\label{JMC}
\begin{aligned}
&\J_i^{\textrm{\tiny \rm MC}}(x_0;{u}_i)\=\frac{1}{2}\E \Bigg\{{\T}\Big[\|x_i\-\Gamma_1\Eo x_i
\|_{Q}^{2} \+ \|u_i\|^{2}_{R}\Big]dt\\
&\hspace{3cm} +
\|x_i(T)\-\Gamma_2\Eo x_i(T)\|_{G}^{2}\Bigg\}.\\
\end{aligned}
\end{equation}
Then we propose the MC problem.

\quad\\
(\textbf{MC}): For given initial value $x_0\in\R^n$, each agent $\A_i$ finds a strategy $\u_i$   such that
\[\J^{\textrm{\tiny \rm MC}}_i(x_0;\u_i) = \inf_{u_i \in \U_i^{ol}}\J_i^{\textrm{\tiny \rm MC}}(x_0;u_i),\]
 where $(x_i,u_i)$ is subject to \eqref{E x}.

\quad\\
Similarly we have the optimal strategy  $\u_i$, optimal state process set $\x_i$ and the optimal pair $(\x_i, \u_i)$.

Next, we present some assumptions to the coefficients which will be in force throughout this paper.
\begin{enumerate}
  \item[(\textbf{A1})] (State-coefficients): $A$, $\bar{A}$, $C$, $\bar{C}$, $C_0$, $\bar{C}_0$ $\in$ $L^\infty(0,T;\R^{n\times n})$, $B$, $\bar{B}$, $D$, $\bar{D}$, $D_0$, $\bar{D}_0$ $\in$ $L^\infty(0,T;\R^{n\times m})$.
  \item[(\textbf{A2})] (Weight-coefficients):  $Q$  $\in$ $L^\infty(0,T;\S^{n })$, $R$ $\in$ $L^\infty(0,T;\S^{m})$, $G$ $\in$ $ \S^{n} $, $\Gamma_1$ $\in$ $L^\infty(0,T;\R^{n\times n})$, $\Gamma_2$ $\in$ $\R^{n\times n}$.
\end{enumerate}
\begin{remark}
  The assumptions above are commonly used in LQG  models, and readers
can refer \cite{yong2013linear, graber2016linear,
bensoussan2016linear} for more applications. Under {(\textbf{A1})},
for any given $\bu \in \U_c^{ol}\times\cdots\times \U_c^{ol}$, \eqref{xN} and \eqref{E x} both admit unique
strong solutions % by Proposition 2.6 in
(see e.g., \cite{yong2013linear}). Furthermore, under {(\textbf{A2})},
the cost functionals $\J_i^{\textrm{\tiny MG}}$, $\J^{\textrm{\tiny \rm MT}}$ and $\J_i^{\textrm{\tiny \rm MC}}$ are well-defined for all $\bu\in \U_c^{ol} \times\cdots\times \U_c^{ol}$.
\end{remark}

\section{Relations among (MC), (MG) and (MT)} \label{sec 4}

In this section, we investigate the relations among (MC), (MG) and (MT). As we know, the problems have been studied in the literature. In particular, (MC) can be found in \cite{yong2013linear}, (MG) in \cite{huang2007large}, and (MT) in \cite{HCM2012}. Here, we first state some systems involved in the procedure of decentralized strategies.

\subsection{Existing results}

\subsubsection{(MC)}
For the sake of notation simplicity, we denote $\widehat{Q} = (I-\Gamma_1)^{\top}Q(I-\Gamma_1)$, $\widehat{G} = (I-\Gamma_2)^{\top}G(I-\Gamma_2)$, $\A = A+\bar{A}$, $\mathcal{B} = B+\bar{B}$, $\mathcal{C} = C+\bar{C}$, $\mathcal{D} = D+\bar{D}$, $\mathcal{C}_0 = C_0+\bar{C}_0$, $\mathcal{D}_0 = D_0+\bar{D}_0$, $\mathcal{R}_1  = (D^{\top}P_1  D \+ D^{\top}_0P_1  D_0   \+ R)$, $ \mathcal{R}_2  = (\mD^{\top}P_1  \mD \+ \mD^{\top}_0P_2  \mD _0 \+ R) $.
Introducing two Riccati equations {\rm(RE1)} and {\rm(RE2)}\s
\begin{equation}
\text{{\rm(RE1)}}\left\{\begin{aligned}
 &  \dot{P}_1    \+P_1 A \+ A^{\top}P_1 \+ C^{\top}P_1 C \+ C^{\top}_0P_1 C_0 \- (P_1 B\\
 &\ \ \+ C^{\top}P_1 D\+ C^{\top}_0P_1 D_0)\mathcal{R}_1^{-1}(B^{\top}P_1  \+  D^{\top}P_1 C\\
 &\ \ \+  D^{\top}_0P_1 C_0) \+ Q  = 0,\quad
   P_1 (T) = G,\\
\end{aligned}\right.
\end{equation}
\begin{equation}
\text{{\rm(RE2)}}\left\{\begin{aligned}
&  \dot{P}_2  \+  P_2 \A  \+ \A^{\top} P_2 \+   \C^{\top} P_1 \C \+   \C^{\top}_0 P_2 \C_0\- (P_2 \B \\
& \ \  \+  \C^{\top} P_1 \mD \+  \C_0^{\top} P_2 \mD_0) \mathcal{R}_2^{-1}(\B^{\top}P_2  \+   \mD^{\top}P_1 \C\\
& \ \  \+   \mD^{\top}_0P_2 \C_0)\+ \widehat{Q} = 0,  \quad P_2 (T) = \widehat{G},\\
\end{aligned}\right.
\end{equation}\n
and Assumption (\textbf{A3})
\begin{enumerate}
  \item[(\textbf{A3})] {\rm(RE1)},  {\rm(RE2)} admit  solutions $P_1,P_2\in C([0,T];\S^{n})$ such that
      \[\mathcal{R}_1(t), \ \mathcal{R}_2(t)\gg 0, \text{ a.e. } t\in[0,T].\]
\end{enumerate}

Under {\rm(\textbf{A1})-(\textbf{A3})}, the functional $u_i\mapsto \J_i^{\textrm{\rm MC}}(x_0;u_i)$ is uniformly convex, and hence {\rm(\textbf{MC})} admits a unique optimal pair $(\x_i ,\u_i )$ on $\U_i$. The following Hamiltonian system (\hmc)
  \begin{equation*}
(\hmc)\left\{\begin{aligned}
 & d\x_i  \= (A\x_i  \+ \bar{A}\Eo \x_i  \+ B \u_i  \+ \bar{B}\Eo \u_i  )dt\\
 &\hspace{10mm}  \+ (C\x_i  \+ \bar{C}\Eo \x_i  \+ D \u_i  \+ \bar{D}\Eo \u_i )dW_i \\
  &\hspace{10mm}  \+ (C_0\x_i  \+ \bar{C}_0\Eo \x_i  \+ D_0 \u_i  \+ \bar{D}_0\Eo \u_i )dW_0,\\
 & d k_i \= -\Big( Q\x_i \- (Q\Gamma_1 +  \Gamma_1^{\top}Q\-\Gamma_1^{\top}Q\Gamma_1)\Eo \x_i  \+ A^{\top}k_i\\ &\hspace{1.3cm}\+\bar{A}^{\top}\Eo k_i \+ C^{\top}\zeta_i \+ \bar{C}^{\top}\Eo \zeta_i \+ C^{\top}_0\zeta_i^0 \\
 &\hspace{1.3cm}\+ \bar{C}^{\top}_0\Eo \zeta_i^0 \Big) dt  \+ \zeta_i dW_i \+ \zeta_i^0 dW_0 , \\
 &\x_i (0) \= x_0, \\
 &k_i(T) \= G\x_i (T)\- (G\Gamma_2 +  \Gamma_2^{\top}G\-\Gamma_2^{\top}G\Gamma_2)\Eo \x_i (T),\\
 &  B^{\top}k_i \+ \bar{B}^{\top}\Eo k_i \+ D^{\top}\zeta_i + \bar{D}^{\top}\Eo \zeta_i \+ D^{\top}_0\zeta_i^0\\
 &\  + \bar{D}^{\top}_0\Eo \zeta_i^0 \+ R\u \= 0
\end{aligned}\right.
\end{equation*}
admits a unique adapted solution $(\x_i ,k_i,\zeta_i,\zeta_i^0)$. The optimal control and cost functional are given as
  $\u_i \= - \mathcal{R}^{-1}_1(B^{\top}P_1 \+  D^{\top}P_1C \+   D^{\top}_0P_1C_0)(\x_i  \- \Eo \x_i ) - \mathcal{R}^{-1}_2(\B^{\top}P_2 \+   \mD^{\top}P_1\C \+   \mD^{\top}_0P_2\C_0)\Eo \x_i,$
   and $
     \J_i^{\textrm{\tiny\rm MC}}(x_0;\u_i) =  \frac{1}{2}\left\langle P_2(0)x_0,x_0\right\rangle$.

\subsubsection{(MG)}
The limiting problem is given as\\

(\textbf{MG})$^*$:For given initial value $x_0\in\R^n$, each agent $\A_i$ finds a strategy $\breve{u}_i\in \U_i^{ol}$ such that
      \begin{equation}
\begin{aligned}
J_i^{\textrm{\tiny MG}}(x_0;\breve{u}_i) = \inf_{u_i \in \U_i^{ol}}J_i^{\textrm{\tiny MG}}(x_0;u_i),
\end{aligned}
\end{equation}
      where
\begin{equation}
\begin{aligned}
& J_i^{\textrm{\tiny MG}}(x_0;{u}_i)\=\frac{1}{2}\E \Big\{{\T}\Big[\|x_i\-\Gamma_1\bar{m}
\|_{Q}^{2}\\
&\hspace{1cm} \+ \|u_i\|^{2}_{R}\Big]dt+
\|x_i(T)\-\Gamma_2\bar{m}(T)\|_{G}^{2}\Big\}
\end{aligned}
\end{equation}
subject to
\begin{equation}\nonumber
\left\{\begin{aligned}
 & dx_i = (Ax_i \+ \bar{A}\bar{m} \+ Bu_i \+ \bar{B}\bar{w} )dt\\
 &\hspace{1cm} \+ (Cx_i \+
 \bar{C}\bar{m} \+ Du_i \+ \bar{D}\bar{w})dW_i, \\
 &\hspace{1cm}  \+ (C_0x_i \+
 \bar{C}_0\bar{m} \+ D_0u_i \+ \bar{D}_0\bar{w})dW_0, \\
 & x_i(0) = x_0.
 \end{aligned}\right.
\end{equation}
Here, $\bar{m}$, $\bar{w}$ are undetermined stochastic processes at this moment, thus
they should be treated as some exogenous  term. The CC system is
      \begin{equation}\label{mfg cc}
      \begin{aligned}
      \bar{m} = \Eo  \breve x_i(\cdot;x_0,\bar{m},\bar{w}),\quad \bar{w} = \Eo  \breve u_i(\cdot;x_0,\bar{m},\bar{w}).
      \end{aligned}
      \end{equation}

The following Hamiltonian system for {\rm(\textbf{MG})$^*$}
\s
  \begin{equation}
 (\hmg)\left\{\begin{aligned}
 &d\breve{x}_i  = (A\breve{x}_i  \+ \bar{A}\bar{m} \+ B\breve{u}_i  \+ \bar{B}\bar{w} )dt\\
 &\hspace{1cm} \+ (C\breve{x}_i  \+ \bar{C}\bar{m} \+ D\breve{u}_i  \+ \bar{D}\bar{w})dW_i\\
 &\hspace{1cm} \+ (C_0\breve{x}_i  \+ \bar{C}_0\bar{m} \+ D_0\breve{u}_i  \+ \bar{D}_0\bar{w})dW_0, \\
 & d l_i = -\left( A^{\top}l_i \+ C^{\top}\varsigma_i \+ C^{\top}_0\varsigma_i^0 + Q\breve{x}_i  - Q\Gamma_1\bar{m}  \right) dt\\
 &\hspace{1cm} \+ \varsigma_i dW_i \+   \varsigma_i^0 dW_0 , \\
 &\breve{x}_i(0) = x_0,\quad l_i(T) = G\breve{x}_i (T)\- G\Gamma_2\bar{m}(T),\\
 & B^{\top}l_i   \+ D^{\top}\varsigma_i \+ D^{\top}_0\varsigma_i^0  \+ R\breve{u}_i \= 0
\end{aligned}\right.
\end{equation}\n
  admits a unique adapted solution $(\breve{x}_i ,l_i,\varsigma_i,\varsigma_i^0)$ where $(\breve{x}_i ,\breve{u}_i)$ is the optimal pair of {\rm(\textbf{MG})$^*$}.
By plugging \eqref{mfg cc} into (\hmg), one has \s
\begin{equation}
 (\ccmg)\left\{\begin{aligned}
 &d\breve{x}_i  = (A\breve{x}_i  \+ \bar{A}\Eo \breve{x}_i \+ B\breve{u}_i  \+ \bar{B}\Eo\breve{u}_i )dt\\
 &\hspace{1cm} \+ (C\breve{x}_i  \+ \bar{C}\Eo \breve{x}_i \+ D\breve{u}_i  \+ \bar{D}\Eo\breve{u}_i)dW_i\\
 &\hspace{1cm} \+ (C_0\breve{x}_i  \+ \bar{C}_0\Eo\breve{x}_i \+ D_0\breve{u}_i  \+ \bar{D}_0\Eo\breve{u}_i)dW_0, \\
 & d l_i = -\left( A^{\top}l_i \+ C^{\top}\varsigma_i \+ C^{\top}_0\varsigma_i^0 + Q\breve{x}_i  - Q\Gamma_1\Eo\breve{x}_i  \right) dt\\
 &\hspace{1cm} \+ \varsigma_i dW_i \+   \varsigma_i^0 dW_0 , \\
 &\breve{x}_i(0) = x_0,\quad l_i(T) = G\breve{x}_i (T)\- G\Gamma_2\Eo\breve{x}_i (T),\\
 & B^{\top}l_i   \+ D^{\top}\varsigma_i \+ D^{\top}_0\varsigma_i^0  \+ R\breve{u}_i\= 0. \\
\end{aligned}\right.
\end{equation}\n

Introduce an asymmetric Riccati equation
\begin{equation}\label{121}
\text{{\rm(RE3)}}\left\{\begin{aligned}
 &  \dot{P_3} \+  P_3\A\+ A^{\top}P_3 \+ C^{\top}P_1\C  \+ C^{\top}_0P_3\C_0\\
 &   \- (P_3\B\+C^{\top}P_1\mD\+C^{\top}_0P_3\mD_0)\mathcal{R}_3^{-1}\\
 &\times(B^{\top}P_3 \+   D^{\top}P_1\C\+   D^{\top}_0P_3\C_0)     \+ (Q\-Q\Gamma_1)  \= 0 , \\
 &P_3(T) = G(I\- \Gamma_2),
\end{aligned}\right.
\end{equation}
where $\mathcal{R}_3 = (D^{\top}P_1\mD  \+ D^{\top}_0P_3\mD_0  \+ R)$. Under {\rm(\textbf{A1})-(\textbf{A2})}, if Riccati equations {\rm(RE1)}, {\rm(RE3)} admit solutions $P_1\in C([0,T],\S^n)$, $P_3\in C([0,T],\R^n)$ such that $\mathcal{R}_1,\ \mathcal{R}_3 \gg 0,\text{ a.e. }t\in[0,T]$, we get the closed-loop system\s
  \begin{equation}\label{126_2}
  \left\{\begin{aligned}
  &d\breve{x}_i  = \Big[\left(A \- B  \mathcal{R}_1(B^{\top}P_1 \+  D^{\top}P_1C \+  D^{\top}_0P_1C_0)\right)(\breve{x}_i  \- \Eo \breve{x}_i )\\
  &\+ \left(\A    \- \B\mathcal{R}_3(B^{\top}P_3 \+   D^{\top}P_1\C \+   D^{\top}_0P_3\C_0)\right)\Eo \breve{x}_i\Big]  dt\\
  & \+ \Big[\left(C \- D  \mathcal{R}_1(B^{\top}P_1 \+  D^{\top}P_1C \+  D^{\top}_0P_1C_0)\right)(\breve{x}_i - \Eo\breve{x}_i)\\
  &\+ \left(\C    \- \mD\mathcal{R}_3(B^{\top}P_3 \+   D^{\top}P_1\C \+   D^{\top}_0P_3\C_0)\right)\Eo \breve{x}_i  \Big]   dW_i \\
  & \+ \Big[\left(C_0 \- D_0  \mathcal{R}_1(B^{\top}P_1 \+  D^{\top}P_1C \+  D^{\top}_0P_1C_0)\right)(\breve{x}_i - \Eo\breve{x}_i)\\
  &\+ \left(\C_0    \- \mD_0\mathcal{R}_3(B^{\top}P_3 \+   D^{\top}P_1\C \+   D^{\top}_0P_3\C_0)\right)\Eo \breve{x}_i  \Big]   dW_0, \\
  &\breve{x}_i(0) = x_0,
  \end{aligned}\right.
    \end{equation}\n
$\breve{u}_i = - \mathcal{R}_1(B^{\top}P_1 \+  D^{\top}P_1C \+  D^{\top}_0P_1C_0)(\breve{x}_i  \- \Eo \breve{x}_i )\- \mathcal{R}_3(B^{\top}P_3 \+   D^{\top}P_1\C \+   D^{\top}_0P_3\C_0)\Eo \breve{x}_i,$
  and $\bar{m}\=\Eo \breve{x}_i$, $\bar{w} \=  \- \mathcal{R}_3(B^{\top}P_3 \+   D^{\top}P_1\C \+   D^{\top}_0P_1\C_0)\bar{m} $.

\subsubsection{(MT)}

Based on person-by-person optimality method, one introduces the following auxiliary problem for $\A_i$:

\quad\\
(\textbf{MT})$^*$: For given initial value $x_0$, each agent $\A_i$ finds a strategy $\check{u}_i \in \mathcal{U}_i^{ol}$ such that $J_i^{\textrm{\tiny \rm MT}}(x_0;\check{u}_i) = \inf_{u_i \in \mathcal{U}_i^{ol}}J_i^{\textrm{\tiny \rm MT}}(x_0;u_i)$, where
\begin{align}
\notag& {J}_i^{\textrm{\tiny \rm MT}}(x_0;u_i)\=\frac{1}{2}\E \bigg\{\int_0^{T} \Big[\|x_i\|^{2}_{Q} \+ 2\left\langle M_1, x_i\right\rangle \+ 2\left\langle M_2, u_i\right\rangle\\
&\hspace{2cm} \+ \|u_i\|^{2}_{R}\Big]dt+ \|x_i(T)\|^{2}_{G} + 2\left\langle M_3, x_i(T)\right\rangle\bigg\}, \\
& \text{s.t. }\left\{\begin{aligned}
 & dx_i = (Ax_i \+ \bar{A}\bar{m} \+ Bu_i \+ \bar{B}\bar{w} )dt\\
 &\hspace{1cm} \+ (Cx_i \+
 \bar{C}\bar{m} \+ Du_i \+ \bar{D}\bar{w})dW_i \\
 &\hspace{1cm}  \+ (C_0x_i \+
 \bar{C}_0\bar{m} \+ D_0u_i \+ \bar{D}_0\bar{w})dW_0, \\
 & x_i(0) = x_0,
 \end{aligned}\right.
\end{align}
$M_1$, $M_2$ and $M_3$ are given by
\begin{equation*}
\left\{\begin{aligned}
 &  M_1 = \left(\Gamma^{\top}_1Q\Gamma_1  \- \Gamma^{\top}_1Q \- Q\Gamma_1\right) \bar{m} \+ \bar{A}^{\top}\Eo p_j^{1}\\
 &\hspace{1.2cm}\+ \bar{C}^{\top}\Eo q_j^{1} \+ \bar{C}^{\top}_0\Eo q_j^{1,0} \+ \bar{C}^{\top}_0q_0^2 \+ \bar{A}^{\top}p^{2}  , \\
 &  M_2 =  \bar{B}^{\top}\Eo p_j^{1} \+ \bar{D}^{\top}\Eo q_j^{1} \+ \bar{D}^{\top}_0\Eo q_j^{1,0} \+ \bar{D}^{\top}_0q_0^2 \+ \bar{B}^{\top}p^{2},\\  & M_3 = (\Gamma^{\top}_2G\Gamma_2   \- \Gamma^{\top}_2G \- G\Gamma_2)\bar{m}(T).
\end{aligned}\right.
\end{equation*}
Besides, the Hamiltonian system of {\rm(\textbf{MT})$^*$} is
  \begin{equation}
 (\hmt)\left\{\begin{aligned}
 & d\check{x}_i = (A\check{x}_i \+ B\check{u}_i \+ \bar{A}\bar{m} \+ \bar{B}\bar{w})dt\\
 &\hspace{0.5cm} \+ (C\check{x}_i \+ D\check{u}_i \+ \bar{C}\bar{m} \+ \bar{D}\bar{w})dW_i\\
 &\hspace{0.5cm}\+ (C_0\check{x}_i \+ D_0\check{u}_i \+ \bar{C}_0\bar{m} \+ \bar{D}_0\bar{w})dW_0,\\
 & d p_i = \Big[\- Q\check{x}_i \+ Q\Gamma_1 \bar{m} \- A^{\top}p_i \- C^{\top}q_i \- C^{\top}_0q_i^0\\
 &\hspace{0.5cm} \+
 \Gamma^{\top}_1Q(I - \Gamma) \bar{m} \- \bar{A}^{\top}p^2 - \bar{C}^{\top}\Eo q^1_j \- \bar{A}^{\top}\Eo p^1_j \\
 &\hspace{0.5cm}
 \- \bar{C}^{\top}_0\Eo q_j^{1,0}  \+ \bar{C}^{\top}_0q_0^2\Big] dt \+ q_i dW_i  \+ q_i^0 dW_0, \\
 & \check{x}_i(0) = x_0,\quad p_i(T) = G\check{x}_i(T) + M_3,\\
 &  R\check{u}_i \+ B^{\top}p_i \+ D^{\top}q_i  \+ D^{\top}_0q_i^0 \+
\bar{B}^{\top}\Eo p_j^{1} \+ \bar{D}^{\top}\Eo q_j^{1} \\
&\ \ \+ \bar{D}^{\top}_0\Eo q_j^{1,0} \+
\bar{B}^{\top}p^{2} \+ \bar{D}^{\top}_0q_0^2   \= 0.
 \end{aligned}\right.
\end{equation}
Furthermore, one derives CC system\s
\begin{equation}\label{56}
 (\ccmt)\left\{ \begin{aligned}
 & d\check{x}_i = (A\check{x}_i \+ B\check{u}_i \+ \bar{A}\bar{m} \+ \bar{B}\bar{w})dt\\
  &\hspace{1cm}\+ (C\check{x}_i \+ D\check{u}_i \+ \bar{C}\bar{m} \+ \bar{D}\bar{w})dW_i\\
 &\hspace{1cm}\+ (C_0\check{x}_i \+ D_0\check{u}_i \+ \bar{C}_0\bar{m} \+ \bar{D}_0\bar{w})dW_0,\\
 & d p_i \= \Big[- Q\left(\check{x}_i \- \Eo\check{x}_i\right)\- \hat{Q} \Eo\check{x}_i \- A^{\top}p_i \- C^{\top}q_i \\
 &\hspace{1cm}\- C^{\top}_0q_i^0 \- \bar{A}^{\top}p^2 \- \bar{C}^{\top}\Eo q^1_i \- \bar{C}^{\top}_0\Eo q^{1,0}_i\\
 &\hspace{1cm}\- \bar{A}^{\top}\Eo p^1_i \- \bar{C}^{\top}_0q_0^2 \Big] dt \+ q_i dW_i \+ q_i^0 dW_0, \\
 & dp_i^1 \= - \left(Q\check{x}_i \+ A^{\top}p_i^{1} \+ C^{\top}q_i^{1} \+ C^{\top}_0q_i^{1,0}\right) dt\\
 &\hspace{1cm} \+ q_i^{1} dW_i \+ q_i^{1,0} dW_0, \\
 & dp^{2} \= -\Big[(\Gamma^{\top}_1Q\Gamma_1  \- \Gamma^{\top}_1Q\- Q\Gamma_1) \Eo\check{x}_i \+ \bar{A}^{\top}\Eo p_i^{1}\\
 &\hspace{1cm}\+ \bar{C}^{\top}\Eo q_i^{1} \+ \bar{C}^{\top}_0\Eo q_i^{1,0} \+ \bar{A}^{\top}p^{2} \+ A^{\top}p^{2}\\
 &\hspace{1cm}\+ {\C}^{\top}_0q_0^2 \Big] dt \+ q^{2}_0dW_0, \\
 & \check{x}_i(0) \= x_0, \quad p_i^1(T) \= G\check{x} _{i}(T),\\
 &p_i(T) \= G\check{x}_i(T) + (\Gamma^{\top}_2G\Gamma_2   \- \Gamma^{\top}_2G \- G\Gamma_2)\Eo\check{x}_i(T),\\
 & p^{2}(T) \= (\Gamma^{\top}_2G\Gamma_2   \- \Gamma^{\top}_2G \- G\Gamma_2)\Eo\check{x}_i(T),\\
 &    R\check{u}_i \+ B^{\top}p_i \+ D^{\top}q_i \+ D^{\top}_0q_i^0 \+
\bar{B}^{\top}\Eo p_i^{1} \+ \bar{D}^{\top}\Eo q_i^{1}\\
& \+ \bar{D}^{\top}_0\Eo q_i^{1,0} \+
\bar{B}^{\top}p^{2}\+ \bar{D}^{\top}_0q_0^2    \= 0.\\
 \end{aligned}
 \right.
\end{equation}\n

In addition, under {\rm(\textbf{A1})-(\textbf{A3})}, we obtain the closed-loop system\s
  \begin{equation}
  \left\{\begin{aligned}
  &d\check{x}_i  = \Big(\left[A    \- B  \mathcal{R}_1^{-1}(B^{\top}P_1 \+  D^{\top}P_1C \+  D^{\top}_0P_1C_0)\right](\check{x}_i  \- \Eo \check{x}_i )\\
  &\hspace{5mm} \+ \left[\A  \- \B \mathcal{R}_2^{-1}(\B^{\top}P_2 \+   \mD^{\top}P_1\C \+   \mD^{\top}_0P_2\C_0)\right]\Eo \check{x}_i  \Big)dt\\
  & \hspace{5mm}  \+ \Big(\left[C    \- D \mathcal{R}_1^{-1}(B^{\top}P_1 \+  D^{\top}P_1C\+  D^{\top}_0P_1C_0)\right](\check{x}_i  \- \Eo \check{x}_i )\\
  &\hspace{5mm} \+ \left[\C  \- \mD \mathcal{R}_2^{-1}(\B^{\top}P_2 \+   \mD^{\top}P_1\C \+   \mD^{\top}_0P_2\C_0)\right]\Eo \check{x}_i \Big)dW_i \\
  & \hspace{5mm}  \+ \Big(\left[C_0    \- D_0 \mathcal{R}_1^{-1}(B^{\top}P_1 \+  D^{\top}P_1C\+  D^{\top}_0P_1C_0)\right](\check{x}_i  \- \Eo \check{x}_i )\\
  &\hspace{5mm} \+ \left[\C_0  \- \mD_0 \mathcal{R}_2^{-1}(\B^{\top}P_2 \+   \mD^{\top}P_1\C \+   \mD^{\top}_0P_2\C_0)\right]\Eo \check{x}_i \Big)dW_0, \\
  &\check{x}_i(0) = x_0,
  \end{aligned}\right.
    \end{equation}\n
and $\check{u}_i \= - \mathcal{R}^{-1}_1(B^{\top}P_1 \+  D^{\top}P_1C \+   D^{\top}_0P_1C_0)(\check{x}_i  \- \Eo \check{x}_i )- \mathcal{R}^{-1}_2(\B^{\top}P_2 \+   \mD^{\top}P_1\C \+   \mD^{\top}_0P_2\C_0)\Eo \check{x}_i.$

\subsection{Relations of the decentralized strategies}
\begin{lemma}\label{lem4.1}
 Under {\rm(\textbf{A1})-(\textbf{A3})}, if $\bar{A} = \bar{B} = \bar{C}  = \bar{D} = \bar{C}_0  = \bar{D}_0  = 0$, $\Gamma_1^{\top}Q\-\Gamma_1^{\top}Q\Gamma_1= 0$, $\Gamma_2^{\top}G\-\Gamma_2^{\top}G\Gamma_2= 0$ (e.g., $\Gamma_1 \= \Gamma_2 \= I$ or $\Gamma_1 \= \Gamma_2 \= 0$), then {\rm(\ccmg)} and {\rm(\hmc)} are identical.
\end{lemma}
\begin{proof}
  The result of Lemma \ref{lem4.1} can be obtained  directly by comparing the coefficients of {\rm(\ccmg)} and {\rm(\hmc)}.
\end{proof}
By Lemma \ref{lem4.1} we can obtain the relation between (\textbf{MC}) and (\textbf{MG}) as follows: for each agent $\A_i$, {\rm(\textbf{MC})} admits a unique open-loop optimal strategy denoted by $\u_i \in \U_i^{ol}$, and {\rm(\textbf{MG})} admits a unique open-loop mean-field strategy denoted by $\breve{u}_i  \in \U_i^{ol}$. Moreover in this case $\u_i  = \breve{u}_i $.
%\begin{theorem}\label{them4.1}
%  Let {\rm(\textbf{A1})-(\textbf{A3})} hold and $\bar{A} = \bar{B} = \bar{C}  = \bar{D} = \bar{C}_0  = \bar{D}_0  = 0$, $\Gamma_1^{\top}Q\-\Gamma_1^{\top}Q\Gamma_1= 0$, $\Gamma_2^{\top}G\-\Gamma_2^{\top}G\Gamma_2= 0$ (e.g., $\Gamma_1 \= \Gamma_2 \= I$ or $\Gamma_1 \= \Gamma_2 \= 0$). Then for each agent $\A_i$, {\rm(\textbf{MC})} admits a unique open-loop optimal strategy denoted by $\u_i \in \U_i^{ol}$, and {\rm(\textbf{MG})} admits a unique open-loop mean-field strategy denoted by $\check{u}_i  \in \U_i^{ol}$. Moreover in this case $\u_i  = \check{u}_i $.
%\end{theorem}
%\begin{proof}
%  Under (\textbf{A1})-(\textbf{A3}), by Proposition \ref{prop 3.2}, (\textbf{MC}) admits a unique optimal strategy $\u_i$ which is
%  determined by (\hmc), and (\hmc) admits a unique adapted solution.   Then what we desire to prove next is the existence and uniqueness of MG strategies. By Lemma \ref{lem4.1}, {\rm(\ccmg)} is also uniquely solvable. Thus, the mean-field terms $\bar{m}$ and $\bar{w}$ can be uniquely determined by {\rm(\ccmg)}. Consequently, by Proposition \ref{prop 4.1} and procedure (MG1)-(MG3), there exists a unique open-loop mean-field strategy $\check{u}_i  $ for $\A_i$, which is determined by (\hmg). Last, the relation $\u_i = \check{u}_i$ follows the equivalence of (\ccmg) and (\hmc).
%\end{proof}

Next, we focus on  (\textbf{MC}) and (\textbf{MT}). By comparing (\ccmt) and (\hmc),  we derive
     \begin{lemma}\label{lemma5.1}
  Under {\rm(\textbf{A1})-(\textbf{A2})}, {\rm(\ccmt)} and {\rm(\hmc)} are identical.
\end{lemma}
\begin{proof}
Firstly, we prove (\ccmt) $ \Longrightarrow$ (\hmc) part.

\quad\\
Let $\big(\check{x}_i,p_i,p_i^1, p^2,q_i,q_i^0,q_i^1,q_i^{1,0}\big)$ be the adapted solution of (\ccmt). The dynamic of $p^2 + p^1_i - p_i$ satisfies that\s
\begin{equation}\label{eq911}
\left\{\begin{aligned}
&d(p^2 \+ p^1_i \- p_i)\\
=& -\Big[ A^{\top}(p^{2} \+ p_i^{1}\-p_i)  \+ C^{\top}(q_i^{1}  \- q_i) \+ C^{\top}_0( q_0^2 \+ q_i^{1,0}\-q_i^0)  \Big]dt  \\
& \+ (q_i^{1}\-q_i) dW_i \+ ( q_0^2 \+ q_i^{1,0}\-q_i^0) dW_0 ,\\
&\hspace{-3mm}(p^2 \+ p^1_i \- p_i)(T) \= 0.
\end{aligned}\right.
\end{equation}\n
Under (\textbf{A1})-(\textbf{A2}), \eqref{eq911} admits a unique solution and we have $p^2 + p^1_i - p_i\equiv 0$ and $q_0^2 \+ q_i^{1,0}\-q_i^0 = q_i^{1}-q_i\equiv 0$. Thus, we obtain $\Eo q_i^{1} \equiv \Eo q_i$ and $\Eo q_i^{1,0} \+ q_0^2 \equiv \Eo q_i^{0}$. Next, consider the dynamic of $p^2 + \Eo p^1_i$ which satisfies
\begin{equation}
\left\{\begin{aligned}
&d(p^2 + \Eo p^1_i)\\
=& \Big[\widehat{Q} \Eo\check{x}_i \- \bar{A}^{\top}(\Eo p_i^{1} \+ p^{2}) \- \bar{C}^{\top}\Eo q_i^{1} \- A^{\top}(p^{2} \+ \Eo p_i^{1}) \\
&  \- C^{\top}\Eo q_i^{1} \- \C^{\top}_0\left( \Eo q_i^{1,0} \+ q_0^2\right)\Big] dt \+ \left( \Eo q_i^{1,0} \+ q_0^2\right) dW_0,\\
&\hspace{-3mm} (p^2 + \Eo p^1_i)(T) = \widehat{G}\Eo\check{x}_i(T),
\end{aligned}\right.
\end{equation}
and the dynamic of $\Eo p_i$ satisfying
\begin{equation}
 \left\{\begin{aligned}
 & d\Eo p_i = \Big[\widehat{Q} \Eo\check{x}_i \- A^{\top}\Eo p_i \- C^{\top}\Eo q_i  \- \bar{A}^{\top}(\Eo p_i^{1} \+ p^{2})\\
 &\hspace{13mm} \- \bar{C}^{\top}\Eo q^1_i\- C^{\top}_0\Eo q_i^0 \- \bar{C}^{\top}_0\left( \Eo q_i^{1,0} \+ q_0^2\right)\Big] dt, \\
 &\Eo p_i(T) = \widehat{G}\Eo\check{x}_i(T).
 \end{aligned}\right.
\end{equation}
By taking difference, we obtain\s
\begin{equation}
\left\{\begin{aligned}
&d(p^2 + \Eo p^1_i - \Eo p_i)\\
=& \Big[  \- A^{\top}(p^{2} \+ \Eo p_i^{1} \- \Eo p_i)\- C^{\top}(\Eo q_i^{1} \- \Eo q_i )  \\
& \- C^{\top}_0\left( \Eo q_i^{1,0} \+ q_0^2\right) \Big] dt\+ \left( \Eo q_i^{1,0} \+ q_0^2\right) dW_0,\\
&\hspace{-3mm}(p^2 + \Eo p^1_i - \Eo p_i)(T) = 0.
\end{aligned}\right.
\end{equation}\n
Noticing $q_i^{1}-q_i\equiv 0$, one gets
\begin{equation}
\left\{\begin{aligned}
&d(p^2 + \Eo p^1_i - \Eo p_i)\\
=& \Big[  \- A^{\top}(p^{2} \+ \Eo p_i^{1} \- \Eo p_i) \- C^{\top}_0\left( \Eo q_i^{1,0} \+ q_0^2\right) \Big] dt \\
& + \left( \Eo q_i^{1,0} \+ q_0^2\right) dW_0,\\
&\hspace{-3mm}(p^2 + \Eo p^1_i - \Eo p_i)(T) = 0
\end{aligned}\right.
\end{equation}
and $p^2 + \Eo p^1_i - \Eo p_i \equiv 0$. By plugging  $\Eo q_i^{1} \equiv \Eo q_i$ and $p^2 + \Eo p^1_i - \Eo p_i \equiv 0$ into (\ccmt), we derive
\begin{equation*}
 \left\{\begin{aligned}
 & d\x_i  \= (A\x_i  \+ \bar{A}\Eo \x_i  \+ B \u_i  \+ \bar{B}\Eo \u_i  )dt\\
 &\hspace{0.8cm} \+ (C\x_i  \+ \bar{C}\Eo \x_i  \+ D \u_i  \+ \bar{D}\Eo \u_i )dW_i \\
 &\hspace{0.8cm} \+ (C_0\x_i  \+ \bar{C}_0\Eo \x_i  \+ D_0 \u_i  \+ \bar{D}_0\Eo \u_i )dW_0,\\
 & d p_i \= -\Big( Q\x_i \- (Q\Gamma_1 +  \Gamma_1^{\top}Q\-\Gamma_1^{\top}Q\Gamma_1)\Eo \x_i  \+ A^{\top}p_i \\
 &\hspace{0.8cm}\+\bar{A}^{\top}\Eo p_i \+ C^{\top}q_i \+ \bar{C}^{\top}\Eo q_i \+ C^{\top}_0q_i^0 \+ \bar{C}^{\top}_0\Eo q_i^0 \Big) dt\\
 &\hspace{0.8cm} \+ q_i dW_i \+ q_i^0 dW_0 , \\
 & \x_i (0) \= x_0,\quad p_i(T) \= G\x_i (T)\- (G\Gamma_2 +  \Gamma_2^{\top}G\\
 &\hspace{0.8cm}\-\Gamma_2^{\top}G\Gamma_2)\Eo \x_i (T),\\
 & B^{\top}p_i \+ \bar{B}^{\top}\Eo p_i \+ D^{\top}q_i + \bar{D}^{\top}\Eo q_i \+ D^{\top}_0q_i^0\\
&\hspace{0.8cm} + \bar{D}^{\top}_0\Eo q_i^0 \+ R\u \= 0,
\end{aligned}\right.
\end{equation*}
which is identical to  (\hmc) and $(\check{x}_i, p_i , q_i )$ is a solution of (\hmc).

\quad\\
Secondly, we prove (\hmc) $\Longrightarrow$ (\ccmt) part.

\quad\\
Under (\textbf{A1})-(\textbf{A2}), for any adapted solution $(\bar{x}_i , k_i , \zeta_i, \zeta_i^0 )$ of (\hmc), the following backward stochastic differential equation (BSDE) admits a unique solution \s
\begin{equation} \label{eq88}
  \left\{\begin{aligned}
 & dk_i^1 \= - \left(Q\bar{x}_i \+ A^{\top}k_i^{1} \+ C^{\top}\zeta_i^{1} \+ C^{\top}_0\zeta_i^{1,0}\right) dt \+ \zeta_i^{1} dW_i \+ \zeta_i^{1,0} dW_0, \\
 &k^1_i(T) = G\x_i (T),
 \end{aligned}
 \right.
\end{equation}\n
and so does the following BSDE
\begin{equation} \label{eq89}
  \left\{ \begin{aligned}
 & dk^{2} \= -\Big[(\Gamma^{\top}_1Q\Gamma_1  \- \Gamma^{\top}_1Q\- Q\Gamma_1) \Eo\bar{x}_i \+ \bar{A}^{\top}\Eo k_i^{1} \\
 &\hspace{10mm}\+ \bar{C}^{\top}\Eo \zeta_i^{1}  \+ \bar{C}^{\top}_0\Eo \zeta_i^{1,0}\+ \bar{A}^{\top}k^{2} \+ A^{\top}k^{2}\\
 &\hspace{10mm} \+ {\C}^{\top}_0\zeta_0^2 \Big] dt\+ \zeta^{2}_0dW_0, \\
 &k^{2}(T) = (\Gamma^{\top}_2G\Gamma_2   \- \Gamma^{\top}_2G \- G\Gamma_2)\Eo\x_i (T).
 \end{aligned}
 \right.
\end{equation}
Then by combining (\hmc), \eqref{eq88} and \eqref{eq89}, we  have
\s\begin{equation} \label{eq90}
 \left\{\begin{aligned}
 & d\x_i  \= (A\x_i  \+ \bar{A}\Eo \x_i  \+ B \u_i  \+ \bar{B}\Eo \u_i  )dt\\
 &\hspace{1cm} \+ (C\x_i  \+ \bar{C}\Eo \x_i  \+ D \u_i  \+ \bar{D}\Eo \u_i )dW_i \\
 &\hspace{1cm} \+ (C_0\x_i  \+ \bar{C}_0\Eo \x_i  \+ D_0 \u_i  \+ \bar{D}_0\Eo \u_i )dW_0,\\
 & d k_i \= -\Big[ Q\x_i \- (Q\Gamma_1 +  \Gamma_1^{\top}Q\-\Gamma_1^{\top}Q\Gamma_1)\Eo \x_i  \+ A^{\top}k_i \\
 &\hspace{1cm}\+\bar{A}^{\top}\Eo k_i \+ C^{\top}\zeta_i \+ \bar{C}^{\top}\Eo \zeta_i \+ C^{\top}_0\zeta_i^0 \+ \bar{C}^{\top}_0\Eo \zeta_i^0 \Big] dt\\
 &\hspace{1cm}\+ \zeta_i dW_i \+ \zeta_i^0 dW_0 , \\
 & dk_i^1 \= - \left(Q\bar{x}_i \+ A^{\top}k_i^{1} \+ C^{\top}\zeta_i^{1} \+ C^{\top}_0\zeta_i^{1,0}\right) dt\\
 &\hspace{1cm}\+ \zeta_i^{1} dW_i \+ \zeta_i^{1,0} dW_0, \\
 & dk^{2} \= -\Big[(\Gamma^{\top}_1Q\Gamma_1  \- \Gamma^{\top}_1Q\- Q\Gamma_1) \Eo\bar{x}_i \+ \bar{A}^{\top}\Eo k_i^{1} \+ \bar{C}^{\top}\Eo \zeta_i^{1}\\
 &\hspace{1cm} \+ \bar{C}^{\top}_0\Eo \zeta_i^{1,0} \+ \bar{A}^{\top}k^{2} \+ A^{\top}k^{2} \+ {\C}^{\top}_0\zeta_0^2 \Big] dt \+ \zeta^{2}_0dW_0, \\
 & \x_i (0) \= x_0,\quad k^1_i(T) = G\x_i (T), \\
 &k_i(T) \= G\x_i (T)\- (G\Gamma_2 +  \Gamma_2^{\top}G\-\Gamma_2^{\top}G\Gamma_2)\Eo \x_i (T),\\
 &k^{2}(T) = (\Gamma^{\top}_2G\Gamma_2   \- \Gamma^{\top}_2G \- G\Gamma_2)\Eo\x_i (T),\\
 & B^{\top}k_i \+ \bar{B}^{\top}\Eo k_i \+ D^{\top}\zeta_i + \bar{D}^{\top}\Eo \zeta_i \+ D^{\top}_0\zeta_i^0\\
 &\hspace{1cm}+ \bar{D}^{\top}_0\Eo \zeta_i^0 \+ R\u \= 0.
\end{aligned}\right.
\end{equation}\n
Comparing \eqref{eq90} with (\ccmt), what remains to prove is that $k^2 \+ \Eo k^1_i \equiv \Eo k_i $ and $\Eo \zeta^{1}_i\!\equiv\!\Eo \zeta_i$. Then we consider the dynamic of $k^2 + k^1_i - k_i$
\begin{equation}\label{eq91}
 \left\{ \begin{aligned}
 &d (k^2 \+ k^1_i \- k_i)\\
 =& \- \Big[  A^{\top}(k^{1}_i \+ k^{2} \- k_i) \+ C^{\top}(\zeta^{1}_i \- \zeta_i)\\
 &\+ \bar{A}^{\top}(\Eo k^{1}_i \+ k^{2} \- \Eo k_i) \+ \bar{C}^{\top}(\Eo \zeta^{1}_i   \- \Eo \zeta_i)\\
 & \+ C_0(\zeta^{2}_0\+\zeta_i^{1,0}\-\zeta_i^0)\+ \bar{C}_0(\zeta^{2}_0\+\Eo\zeta_i^{1,0}\-\Eo\zeta_i^0)\Big] dt\\
 &\+ (\zeta^{1}_i\-\zeta_i) dW_i  \+ (\zeta^{2}_0\+\zeta_i^{1,0}\-\zeta_i^0) dW_0, \\
  =&  \- \Big[  A^{\top}(k^{1}_i \+ k^{2} \- k_i) \+ C^{\top}(\zeta^{1}_i \- \zeta_i) \+ \bar{A}^{\top}(\Eo k^{1}_i \+ k^{2} \\
 &\- \Eo k_i)\+ \bar{C}^{\top}\Eo( \zeta^{1}_i   \-  \zeta_i)\+ C_0(\zeta^{2}_0\+\zeta_i^{1,0}\-\zeta_i^0)\\
 &\+ \bar{C}_0\Eo(\zeta^{2}_0\+\zeta_i^{1,0}\-\zeta_i^0)\Big] dt  \+ (\zeta^{1}_i\-\zeta_i) dW_i \\
 &\+ (\zeta^{2}_0\+\zeta_i^{1,0}\-\zeta_i^0) dW_0, \\
 & \hspace{-3mm} (k^2 \+ k^1 \- k)(T) = 0.
 \end{aligned}
 \right.
\end{equation}
Under (\textbf{A1})-(\textbf{A2}), the MF-BSDE \eqref{eq91} admits a unique adapted solution, and obviously $(k^2 \+ k^1_i \- k_i,  \zeta^{1}_i\-\zeta_i) =  (0, 0)$ is an adapted solution of \eqref{eq91}. Thus, $k^2 \+ k^1_i \equiv k_i $, $\zeta^{1}_i\!\equiv\!\zeta_i$, $k^2 \+ \Eo k^1_i \equiv \Eo k_i $ and $\Eo \zeta^{1}_i\!\equiv\!\Eo \zeta_i$. Plugging such relations into \eqref{eq90}, we get
\s\begin{equation}
 \left\{ \begin{aligned}
 & d\bar{x}_i = (A\bar{x}_i \+ B\bar{u}_i \+ \bar{A}\bar{m} \+ \bar{B}\bar{w})dt
 \+ (C\bar{x}_i \+ D\bar{u}_i \+ \bar{C}\bar{m} \\
 &\hspace{1cm}\+ \bar{D}\bar{w})dW_i\+ (C_0\bar{x}_i \+ D_0\bar{u}_i \+ \bar{C}_0\bar{m} \+ \bar{D}_0\bar{w})dW_0,\\
 & d k_i \= \Big[- Q\left(\bar{x}_i \- \Eo\bar{x}_i\right)\- \hat{Q} \Eo\bar{x}_i \- A^{\top}k_i \- C^{\top}\zeta_i \- C^{\top}_0\zeta_i^0\\
 &\hspace{1cm}\- \bar{A}^{\top}k^2 \- \bar{C}^{\top}\Eo \zeta^1_i \- \bar{C}^{\top}_0\Eo \zeta^{1,0}_i\- \bar{A}^{\top}\Eo k^1_i \- \bar{C}^{\top}_0\zeta_0^2 \Big] dt \\
 &\hspace{1cm}\+ \zeta_i dW_i \+ \zeta_i^0 dW_0, \\
 & dk_i^1 \= - \left(Q\bar{x}_i \+ A^{\top}k_i^{1} \+ C^{\top}\zeta_i^{1} \+ C^{\top}_0\zeta_i^{1,0}\right) dt \+ \zeta_i^{1} dW_i \+ \zeta_i^{1,0} dW_0, \\
 & dk^{2} \= -\Big[(\Gamma^{\top}_1Q\Gamma_1  \- \Gamma^{\top}_1Q\- Q\Gamma_1) \Eo\bar{x}_i \+ \bar{A}^{\top}\Eo k_i^{1} \+ \bar{C}^{\top}\Eo \zeta_i^{1}  \\
 &\hspace{10mm}\+ \bar{C}^{\top}_0\Eo \zeta_i^{1,0} \+ \bar{A}^{\top}k^{2} \+ A^{\top}k^{2} \+ {\C}^{\top}_0\zeta_0^2 \Big] dt \+ \zeta^{2}_0dW_0, \\
 & \bar{x}_i(0) \= x_0,\quad k_i^1(T) \= G\bar{x} _{i}(T), \\
 &k_i(T) \= G\bar{x}_i(T) + (\Gamma^{\top}_2G\Gamma_2   \- \Gamma^{\top}_2G \- G\Gamma_2)\Eo\bar{x}_i(T),\\
 &k^{2}(T) \= (\Gamma^{\top}_2G\Gamma_2   \- \Gamma^{\top}_2G \- G\Gamma_2)\Eo\bar{x}_i(T),\\
 & R\bar{u}_i \+ B^{\top}k_i \+ D^{\top}\zeta_i \+ D^{\top}_0\zeta_i^0 \+
\bar{B}^{\top}\Eo k_i^{1} \+ \bar{D}^{\top}\Eo \zeta_i^{1}\\
& \hspace{1cm}\+ \bar{D}^{\top}_0\Eo \zeta_i^{1,0} \+
\bar{B}^{\top}k^{2}\+ \bar{D}^{\top}_0\zeta_0^2 \= 0,
 \end{aligned}
 \right.
\end{equation}\n
which is identical to (\ccmt) and $\big(\bar{x}_i,k_i,k_i^1,k^2,\zeta_i,\zeta_i^0,\zeta_i^1,\zeta_i^{1,0}\big)$ is a solution of (\ccmt), where $\big(k_i^1,\zeta_i^1,\zeta_i^{1,0}\big)$  and $\left(k_2,\zeta_0^2\right)$ are the unique solutions of \eqref{eq88}  and \eqref{eq89} respectively. This completes the proof.
\end{proof}
By the discussion above, we can conclude the following result of the equivalence between the MC optimal control and mean-field control of (\textbf{MT}).
\begin{theorem}\label{theorem 4.1}
  Under {\rm(\textbf{A1})-(\textbf{A3})}, for each agent $\A_i$, {\rm(\textbf{MC})} admits a unique open-loop optimal strategy   ${\u}_i \in \U_i^{ol}$, and {\rm(\textbf{MT})} admits a unique open-loop mean-field strategy denoted by $\check{u}_i \in \U_i^{ol}$. Moreover $\u_i = \check{u}_i $.
\end{theorem}
We can conclude the relations as the following graph, where the arrows ``$\rightarrow$" means ``imply" or ``derive".
\begin{figure}[H]
\begin{center}
\includegraphics[width=\linewidth]{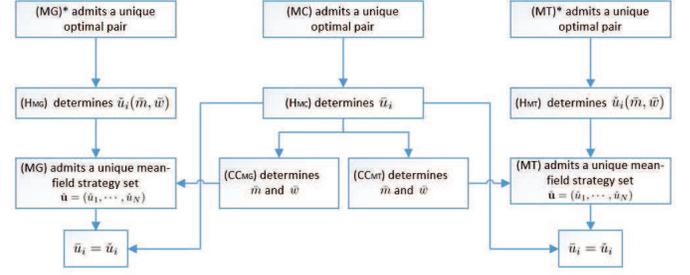}
\caption{Relations among {\rm(\textbf{MG})}, {\rm(\textbf{MC})} and {\rm(\textbf{MT})}}
  \end{center}
\end{figure}
\subsection{Relations of the value functions}
\begin{proposition}\label{prop 4.5}
  Let $\bar{\bu}  := (\bar{u}_1,\cdots,\bar{u}_N)$ be the feedback optimal strategy set of {\rm(\textbf{MC})}  and $\tilde{\bu} := (\tilde{u}_1,\cdots,\tilde{u}_N)$ be the feedback mean-field centralized strategy set of {\rm(\textbf{MT})}. Then
  \begin{equation}
    \|\J_i^{\textrm{\tiny \rm MC}}(x_0;\bar{u}_i) - \J_i^{\textrm{\tiny \rm MT}}(x_0;\tilde{\bu})\| = O\left(\frac{1}{\sqrt{N}}\right), \text{ for any } i\in\I.
  \end{equation}
\end{proposition}
\begin{proof}
Define $\Theta_1:=- \mathcal{R}^{-1}_2(\B^{\top}P_2 \+   \mD^{\top}P_1\C \+   \mD^{\top}_0P_2\C_0),\ \Theta_2:=- \mathcal{R}^{-1}_1(B^{\top}P_1 \+  D^{\top}P_1C \+   D^{\top}_0P_1C_0)$. Then the feedback mean-field strategies  of (\textbf{MT}) and (\textbf{MC}) are
\begin{equation}
  \begin{aligned}
  &\tilde{u}_i =\Theta_1\Eo \tilde{x}_i \+ \Theta_2(\tilde{x}_i  - \Eo \tilde{x}_i ),
  \end{aligned}
  \end{equation}
   and
\begin{equation}
  \begin{aligned}
   \u_i  =&\Theta_1\Eo \bar{x}_i \+ \Theta_2(\bar{x}_i  - \Eo \bar{x}_i ),
  \end{aligned}
  \end{equation}
  where $\tilde{x}_i$ is the realized states of (\textbf{MT}).
  %satisfying
%\begin{equation}\nonumber
%  \left\{\begin{aligned}
%  &d\x_i  = \Big\{\Big[A    \- B  \mathcal{R}_1^{-1}(B^{\top}P_1 \+  D^{\top}P_1C \+  D^{\top}_0P_1C_0)\Big](\x_i  \- \Eo \x_i )\\
%  &\hspace{5mm}\+ \left[\A  \- \B \mathcal{R}_2^{-1}(\B^{\top}P_2 \+   \mD^{\top}P_1\C \+   \mD^{\top}_0P_2\C_0)\right]\Eo \x_i  \Big\}dt\\
%  &\hspace{5mm} \+ \Big\{\left[C    \- D \mathcal{R}_1^{-1}(B^{\top}P_1 \+  D^{\top}P_1C\+  D^{\top}_0P_1C_0)\right](\x_i  \- \Eo \x_i )\\
%  &\hspace{5mm}\+ \left[\C  \- \mD \mathcal{R}_2^{-1}(\B^{\top}P_2 \+   \mD^{\top}P_1\C \+   \mD^{\top}_0P_2\C_0)\right]\Eo \x_i \Big\}dW_i \\
%  &\hspace{5mm} \+ \Big\{\left[C_0    \- D_0 \mathcal{R}_1^{-1}(B^{\!\top\!}P_1 \+  D^{\!\top\!}P_1C\+  D^{\!\top\!}_0P_1C_0)\right](\x_i  \- \Eo \x_i )\\
%  &\hspace{5mm}\+ \left[\C_0  \- \mD_0 \mathcal{R}_2^{-1}(\B^{\top}P_2 \+   \mD^{\top}P_1\C \+   \mD^{\top}_0P_2\C_0)\right]\Eo \x_i \Big\}dW_0, \\
%  &\x_i(0) = x_0.
%  \end{aligned}\right.
%    \end{equation}\n

  Firstly, we estimate $\Eo \tilde{x}_i \- \Eo\bar{x}_i$. We derive
  \begin{equation}
  \left\{\begin{aligned}
& d\Eo\tilde{x}_i  =(A\Eo\tilde{x}_i \+ \bar{A}\Eo\tilde{x}^{(N)}  \+ B\Eo\tilde{u}_i \+ \bar{B}\Eo\tilde{u}^{(N)} )dt\\
 &\+ (C_0\Eo\tilde{x}_i \+ \bar{C}_0\Eo\tilde{x}^{(N)}  \+ D_0\Eo\tilde{u}_i \+ \bar{D}_0\Eo\tilde{u}^{(N)} )dW_0 \\
 &\qquad\ = (A\Eo\tilde{x}_i \+ \bar{A}\Eo\tilde{x}_i  \+ B\Theta_1\Eo \tilde{x}_i \+ \bar{B}\Theta_1\Eo \tilde{x}_i )dt\\
  &\+ (C_0\Eo\tilde{x}_i \+ \bar{C}_0\Eo\tilde{x}_i  \+ D_0\Theta_1\Eo \tilde{x}_i \+ \bar{D}_0\Theta_1\Eo \tilde{x}_i )dW_0,  \\
& \Eo\tilde{x}_i(0) = x_0,
 \end{aligned}\right.
  \end{equation}
  and
  \begin{equation}
  \left\{\begin{aligned}
  &d\Eo\x_i    = (A\Eo\bar{x}_i \+ \bar{A}\Eo\bar{x}_i  \+ B\Theta_1\Eo \bar{x}_i \+ \bar{B}\Theta_1\Eo \bar{x}_i )dt\\
  & \+ (C_0\Eo\bar{x}_i \+ \bar{C}_0\Eo\bar{x}_i  \+ D_0\Theta_1\Eo \bar{x}_i \+ \bar{D}_0\Theta_1\Eo \bar{x}_i )dW_0,  \\
  &\Eo\x_i(0) = x_0.
  \end{aligned}\right.
    \end{equation}
 Thus, $\Eo \tilde{x}_i \= \Eo\bar{x}_i$.

 Secondly, we estimate $\tilde{x}^{(N)} \- \Eo\bar{x}_i$. We have
 \begin{equation*}
 \begin{aligned}
 \tilde{u}^{(N)} = \Theta_1\Eo \tilde{x}_i \+ \Theta_2(\tilde{x}^{(N)}  - \Eo \tilde{x}_i ),
 \end{aligned}
 \end{equation*}
 and
 \begin{equation}
  \left\{\begin{aligned}
 & d\tilde{x}^{(N)} = \big[\A\tilde{x}^{(N)}   \+ \B\Theta_1\Eo \tilde{x}_i \+ \B\Theta_2(\tilde{x}^{(N)}  - \Eo \tilde{x}_i )\big]dt\\
  &\hspace{0.2cm}\+ \frac{1}{N}\sum_{i=1}^{N}(C\tilde{x}_i \+ \bar{C}\tilde{x}^{(N)}  \+ D\tilde{u}_i \+ \bar{D}\tilde{u}^{(N)})dW_i, \\
 &\hspace{0.2cm}\+  \big[\C_0\tilde{x}^{(N)}   \+ \mD_0\Theta_1\Eo \tilde{x}_i \+ \mD_0\Theta_2(\tilde{x}^{(N)}  - \Eo \tilde{x}_i ) \big]dW_0, \\
 & \tilde{x}_i(0) = x_0.
 \end{aligned}\right.
  \end{equation}
  Then
    \begin{equation*}
    \left\{\begin{aligned}
    &d\left(\tilde{x}^{(N)} \- \Eo\bar{x}_i \right)\\
    %=& \Big[ \A\left(\tilde{x}^{(N)} \- \Eo\bar{x}_i \right)  \+ \B \Theta_1\left(\tilde{x}^{(N)} \- \Eo\bar{x}_i \right)\\
%    & + \B\Theta_2(\tilde{x}^{(N)}  \- \Eo \tilde{x}_i )\Big] dt \+ \frac{1}{N}\sum_{i=1}^{N}(C\tilde{x}_i \+ \bar{C}\tilde{x}^{(N)}  \+ D\tilde{u}_i \\
%    &  \+ \bar{D}\tilde{u}^{(N)})dW_i \+ \Big[\C_0\left(\tilde{x}^{(N)} \- \Eo\bar{x}_i \right)   \\
%    &\+ \mD_0 \Theta_1\left(\tilde{x}^{(N)} \- \Eo\bar{x}_i \right) \+ \mD_0\Theta_2(\tilde{x}^{(N)}  \- \Eo \tilde{x}_i )\Big] dW_0, \\
    %=& \Big[\A\left(\tilde{x}^{(N)} \- \Eo\bar{x}_i \right)  \+ \B \Theta_1\left(\tilde{x}^{(N)} \- \Eo\bar{x}_i \right)\\
%    &\+ \B\Theta_2(\tilde{x}^{(N)} \- \Eo\bar{x}_i\+ \Eo\bar{x}_i  - \Eo \tilde{x}_i )\Big] dt\\
%    &\+ \frac{1}{N}\sum_{i=1}^{N}(C\tilde{x}_i \+ \bar{C}\tilde{x}^{(N)}  \+ D\tilde{u}_i \+ \bar{D}\tilde{u}^{(N)})dW_i, \\
%    &\+\Big[\C_0\left(\tilde{x}^{(N)} \- \Eo\bar{x}_i \right) \+ \mD_0 \Theta_1\left(\tilde{x}^{(N)} \- \Eo\bar{x}_i \right)\\
%    & \+ \mD_0\Theta_2(\tilde{x}^{(N)} \- \Eo\bar{x}_i\+ \Eo\bar{x}_i  - \Eo \tilde{x}_i )\Big] dW_0\\
    = & \Big[\A\left(\tilde{x}^{(N)} \- \Eo\bar{x}_i \right)  \+ \B \Theta_1\left(\tilde{x}^{(N)} \- \Eo\bar{x}_i \right)\\
    &+ \B\Theta_2(\tilde{x}^{(N)} \- \Eo\bar{x}_i )\Big] dt \+ \frac{1}{N}\sum_{i=1}^{N}(C\tilde{x}_i \+ \bar{C}\tilde{x}^{(N)}  \\
    & \+ D\tilde{u}_i \+ \bar{D}\tilde{u}^{(N)})dW_i \+ \Big[\C_0\left(\tilde{x}^{(N)} \- \Eo\bar{x}_i \right) \\
    & \+ \mD_0 \Theta_1\left(\tilde{x}^{(N)} \- \Eo\bar{x}_i \right) \+ \mD_0\Theta_2(\tilde{x}^{(N)} \- \Eo\bar{x}_i )\Big] dW_0, \\
    \quad\\
    &\hspace{-3mm}\left(\tilde{x}^{(N)} \- \Eo\bar{x}_i \right)(0) = 0
    \end{aligned}\right.
    \end{equation*}
 which implies $\|\tilde{x}^{(N)} \- \E\bar{x}_i\|^2_{L^2} = O\left(\frac{1}{N}\right)$.

    Thirdly, we estimate $\tilde{x}_i \- \bar{x}_i$. Recall that
    \begin{equation*}
    \left\{\begin{aligned}
    & d\tilde{x}_i = \Big[A\tilde{x}_i \+ \bar{A}\tilde{x}^{(N)}  \+ B\Theta_1\Eo\tilde{x}_i \+ B\Theta_2\Eo(\tilde{x}_i \- \Eo\tilde{x}_i)\\
    &\qquad\+ \bar{B}\Theta_1\Eo\tilde{x}_i \+ \bar{B}\Theta_2(\tilde{x}^{(N)}
     \- \Eo\tilde{x}_i) \Big]dt\\
     &\qquad\+ (C\tilde{x}_i \+ \bar{C}\tilde{x}^{(N)}  \+ D\tilde{u}_i \+ \bar{D}\tilde{u}^{(N)})dW_i \\
     &\qquad \+ \Big[C_0\tilde{x}_i \+ \bar{C}_0\tilde{x}^{(N)}  \+ D_0\Theta_1\Eo\tilde{x}_i \+ D_0\Theta_2\Eo(\tilde{x}_i \- \Eo\tilde{x}_i) \\
    &\qquad \+ \bar{D}_0\Theta_1\Eo\tilde{x}_i\+ \bar{D}_0\Theta_2(\tilde{x}^{(N)}
     \- \Eo\tilde{x}_i) \Big]dW_0,\\
    &d\x_i  = \Big[A\x_i \+ \bar{A}\Eo\x_i \+ B\Theta_1\Eo{\x}_i \+ B\Theta_2({\x}_i \- \Eo{\x}_i)\\
      &\qquad \+ \bar{B}\Theta_1\Eo{\x}_i \Big]dt \+ (C\tilde{x}_i \+ \bar{C}\tilde{x}^{(N)}  \+ D\tilde{u}_i \+ \bar{D}\tilde{u}^{(N)})dW_i\\
      &\qquad \+ \Big[C_0\x_i \+ \bar{C}_0\Eo\x_i \+ D_0\Theta_1\Eo{\x}_i \+ D_0\Theta_2({\x}_i \- \Eo{\x}_i)\\
      &\qquad \+ \bar{D}_0\Theta_1\Eo{\x}_i \Big]dW_0,\\
    &\tilde{x}_i(0) = \x_i(0) = x_0.
    \end{aligned}\right.
    \end{equation*}
    Then
    \begin{equation*}
    \left\{\begin{aligned}
      & d\left(\tilde{x}_i\-\x_i\right)\\
      =& \Big[A\left(\tilde{x}_i\-\x_i\right) \+ \bar{A}(\tilde{x}^{(N)} \- \Eo\bar{x}_i) \\
      &\+ B\Theta_1(\Eo\tilde{x}_i - \Eo\bar{x}_i) \+ B\Theta_2(\left(\tilde{x}_i \- \x_i\right) \- \left(\Eo\tilde{x}_i \- \Eo{\x}_i\right))\\
      &\+ \bar{B}\Theta_1\left(\Eo\tilde{x}_i \- \Eo{\x}_i\right) \+ \bar{B}\Theta_2\Eo(\tilde{x}^{(N)}
     \- \Eo\tilde{x}_i) \Big]dt\\
     &\+\Big[C_0\left(\tilde{x}_i\-\x_i\right) \+ \bar{C}_0(\tilde{x}^{(N)} \- \Eo\bar{x}_i)  \+ D_0\Theta_1(\Eo\tilde{x}_i - \Eo\bar{x}_i)\\
     &\+ D_0\Theta_2(\left(\tilde{x}_i \- \x_i\right) \- \left(\Eo\tilde{x}_i \- \Eo{\x}_i\right))\\
      &\+ \bar{D}_0\Theta_1\left(\Eo\tilde{x}_i \- \Eo{\x}_i\right) \+ \bar{D}_0\Theta_2\Eo(\tilde{x}^{(N)}
     \- \Eo\tilde{x}_i) \Big]dW_0\\
     &  \+ (\cdots)dW_i\\
     = & \Big[(A\+B\Theta_2)\left(\tilde{x}_i\-\x_i\right) \+ O\left(\frac{1}{\sqrt{N}}\right)\Big]dt\\
     &\+ \Big[(C\+D\Theta_2)\left(\tilde{x}_i\-\x_i\right) \+ O\left(\frac{1}{\sqrt{N}}\right)\Big]dW_i\\
     &\+ \Big[(C_0\+D_0\Theta_2)\left(\tilde{x}_i\-\x_i\right) \+ O\left(\frac{1}{\sqrt{N}}\right)\Big]dW_0 ,\\
     \quad\\
     &\hspace{-3mm}\left(\tilde{x}_i\-\x_i\right)(0) = 0
    \end{aligned}\right.
    \end{equation*}
which involves $\|\tilde{x}_i\-\x_i\|^2_{L_2} = O\left(\frac{1}{N}\right)$ and
    \begin{equation}
    \|\J_i^{\textrm{\tiny \rm MC}}(x_0;\bar{u}_i) - \J_i^{\textrm{\tiny \rm MT}}(x_0;\tilde{\bu})\| = O\left(\frac{1}{\sqrt{N}}\right), \text{ for any } i\in\I.
  \end{equation}
\end{proof}
 Next, for the value function of {\rm(\textbf{MG})}$^*$ we have
\begin{proposition}
  Let $\bar{\bu}  := (\bar{u}_1,\cdots,\bar{u}_N)$ be the feedback optimal strategy set of {\rm(\textbf{MC})} and $\breve{\bu} := (\breve{u}_1,\cdots,\breve{u}_N)$ be the feedback optimal decentralized strategy set of {\rm(\textbf{MG})$^*$}. Then
  \begin{equation}
     \J_i^{\textrm{\tiny \rm MC}}(x_0;\bar{u}_i) \leq J_i^{\textrm{\tiny \rm MG}}(x_0;\breve{u}_i), \text{ for any } i\in\I.
  \end{equation}
\end{proposition}
\begin{proof}
  By the fixed point analysis of system (\ccmg),  we have $\bar{m} = \Eo \breve{x}_i$ and $\bar{w} = \Eo \breve{u}_i$. Then it holds that
  \begin{equation}
  \left\{\begin{aligned}
 &d\breve{x}_i  = (A\breve{x}_i  \+ \bar{A}\Eo \breve{x}_i \+ B\breve{u}_i  \+ \bar{B}\Eo\breve{u}_i )dt\\
 &\hspace{1cm} \+ (C\breve{x}_i  \+ \bar{C}\Eo \breve{x}_i \+ D\breve{u}_i  \+ \bar{D}\Eo\breve{u}_i)dW_i\\
 &\hspace{1cm} \+ (C_0\breve{x}_i  \+ \bar{C}_0\Eo\breve{x}_i \+ D_0\breve{u}_i  \+ \bar{D}_0\Eo\breve{u}_i)dW_0, \\
 &\breve{x}_i(0) = x_0,
\end{aligned}\right.
\end{equation}
and
\begin{equation*}
\begin{aligned}
 J_i^{\textrm{\tiny MG}}(x_0;\breve{u}_i)=\ &\frac{1}{2}\E \Big\{{\T}\Big[\|\breve x_i\-\Gamma_1\bar{m}
\|_{Q}^{2} \+ \|\breve u_i\|^{2}_{R}\Big]dt \\
&\+\|\breve x_i(T)\-\Gamma_2\bar{m}(T)\|_{G}^{2}\Big\}\\
=\ &\frac{1}{2}\E \Big\{{\T}\Big[\|\breve x_i\-\Gamma_1\Eo \breve{x}_i
\|_{Q}^{2} \+ \|\breve u_i\|^{2}_{R}\Big]dt\\
& +
\|\breve x_i(T)\-\Gamma_2\Eo \breve{x}_i(T)\|_{G}^{2}\Big\}.
\end{aligned}
\end{equation*}
In this case, we see that  $J_i^{\textrm{\tiny MG}}(x_0;\breve{u}_i) = \J_i^{\textrm{\tiny MC}}(x_0;\breve{u}_i)$. Thus, by the optimality of $\u_i$, it follows
\begin{equation}
     \J_i^{\textrm{\tiny \rm MC}}(x_0;\bar{u}_i) \leq J_i^{\textrm{\tiny \rm MG}}(x_0;\breve{u}_i), \text{ for any } i\in\I.
  \end{equation}
\end{proof}
By the $\varepsilon$-optimality, we also have
\begin{corollary}
  Let $\bar{\bu}  := (\bar{u}_1,\cdots,\bar{u}_N)$ be the feedback optimal strategy set of {\rm(\textbf{MC})}   and $\tilde{\bu} := (\tilde{u}_1,\cdots,\tilde{u}_N)$ be the feedback mean-field centralized strategy set of {\rm(\textbf{MG})}. Then  it holds that
  \begin{equation}
  \begin{aligned}
     &\J_i^{\textrm{\tiny \rm MC}}(x_0;\bar{u}_i)\leq J_i^{\textrm{\tiny \rm MG}}(x_0;\tilde{u}_i) \leq \J_i^{\textrm{\tiny \rm MG}}(x_0;\tilde{\bu})\+O\left(\frac{1}{\sqrt{N}}\right),\\
  \end{aligned}
  \end{equation}
for any $i\in\I$.
\end{corollary}
Thus, combined with Proposition \ref{prop 4.5}, we obtain
\begin{proposition}\label{prop4.8}
  Let $\bar{\bu}  := (\bar{u}_1,\cdots,\bar{u}_N)$ be the feedback optimal strategy set of {\rm(\textbf{MC})}; $\tilde{\bu}^{\text{\rm \tiny MG}} := (\tilde{u}_1^{\text{\rm \tiny MG}},\cdots,\tilde{u}_N^{\text{\rm \tiny MG}})$ and  $\tilde{\bu}^{\text{\rm \tiny MT}} := (\tilde{u}_1^{\text{\rm \tiny MT}},\cdots,\tilde{u}_N^{\text{\rm \tiny MT}})$ be feedback mean-field centralized strategy sets of {\rm(\textbf{MG})} and {\rm(\textbf{MT})}, respectively. Then   it holds that
  \begin{equation}\label{4.22}
  \begin{aligned}
    &\J_i^{\textrm{\tiny \rm MT}}(x_0;\tilde{\bu}^{\text{\rm \tiny MT}})\-O\left(\frac{1}{\sqrt{N}}\right)\leq\J_i^{\textrm{\tiny \rm MC}}(x_0;\bar{u}_i)\leq J_i^{\textrm{\tiny \rm MG}}(x_0;\breve{u}_i)\\
     &\leq \J_i^{\textrm{\tiny \rm MG}}(x_0;\tilde{\bu}^{\text{\rm \tiny MG}})\+O\left(\frac{1}{\sqrt{N}}\right),
  \end{aligned}
  \end{equation}
  for any $i\in\I$.
\end{proposition}
\begin{remark} (i)  We can directly extend the result of Proposition \ref{prop 4.5}-\ref{prop4.8}  from the feedback form decentralized strategy set to open-loop decentralized strategy set.

(ii) Generally, if $\u_i\neq \breve{u}_i$, by the optimality of $\u_i$, it holds that
     $\J_i^{\textrm{\tiny \rm MC}}(x_0;\bar{u}_i) < J_i^{\textrm{\tiny \rm MG}}(x_0;\breve{u}_i)$. In this case, when $N$ is large enough, the relation in \eqref{4.22} can be represented by the following graph
     \begin{figure}[H]
\begin{center}
\includegraphics[width=8cm]{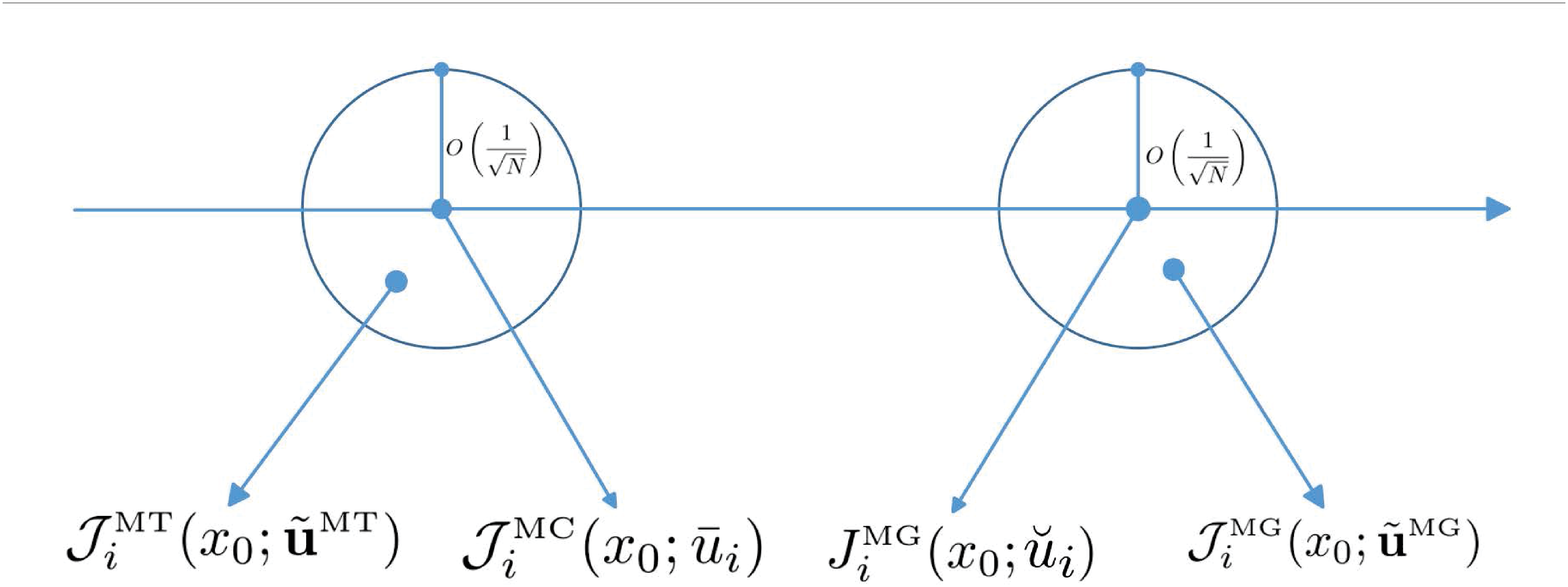}
\caption{Visualization of equation \eqref{4.22} when $\u_i  \neq \breve{u}_i $}\label{fig0}
  \end{center}
\end{figure}
  Moreover, by Lemma \ref{lem4.1}, if $\bar{A} = \bar{B} = \bar{C}  = \bar{D} = \bar{C}_0  = \bar{D}_0  = 0$, $\Gamma_1^{\top}Q\-\Gamma_1^{\top}Q\Gamma_1= 0$, $\Gamma_2^{\top}G\-\Gamma_2^{\top}G\Gamma_2= 0$ (e.g., $\Gamma_1 \= \Gamma_2 \= I$ or $\Gamma_1 \= \Gamma_2 \= 0$), then $\u_i  = \breve{u}_i $. In this case, the relation in \eqref{4.22} can be represented by
  \begin{figure}[H]
\begin{center}
\includegraphics[width=8cm]{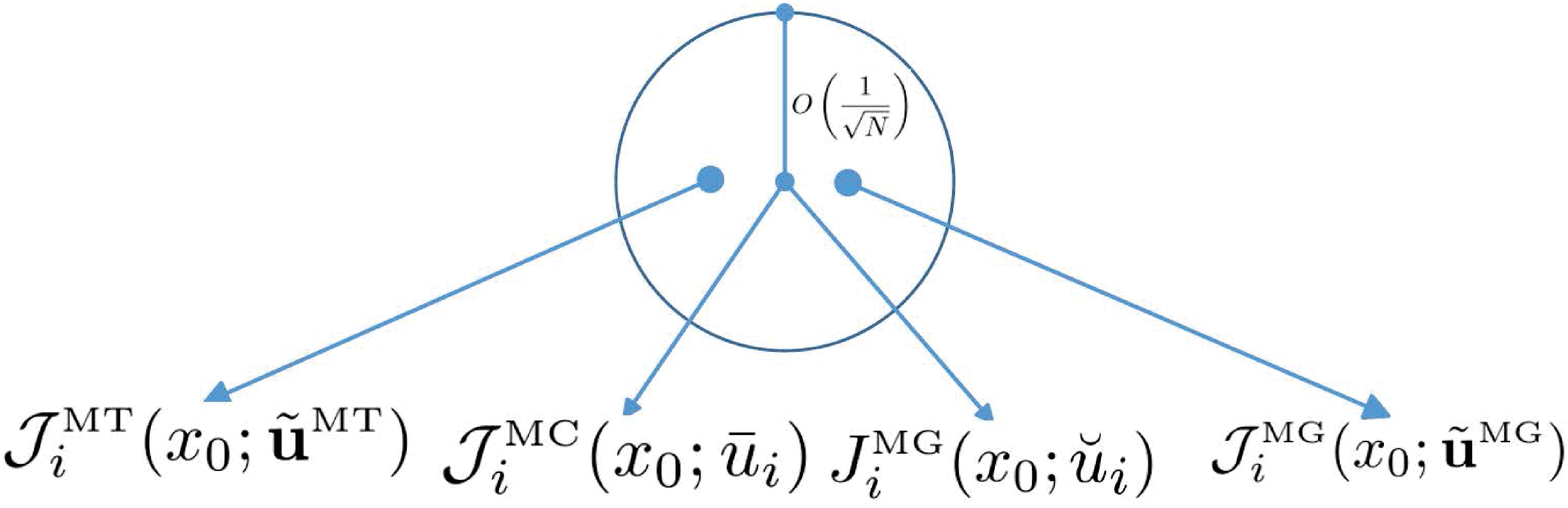}
\caption{Visualization of equation \eqref{4.22} when $\u_i  = \breve{u}_i $ }\label{fig-1}
  \end{center}
\end{figure}

\end{remark}
\subsection{A numerical example}
In this section, we provide a numerical example to illustrate the relations among (\textbf{MC}) (\textbf{MG}) and (\textbf{MT}). We consider the numerical simulation on the time interval $[0,1]$ and the initial state $x_0 = \binom{0.1037}{0.8396}$. For the sake of calculation simplicity, we let $C_0 = \bar C_0 = D_0 = \bar D_0$, and the other coefficients are given as
\textbf{dynamic coefficients:}
\begin{equation*}
\begin{aligned}
&A = \begin{pmatrix}
      0.8903 & 0.6517 \\
      0.1961 & 0.2188
    \end{pmatrix},
\bar{A} = \begin{pmatrix}
      0.7686& 0.8326 \\
      0.0016 & 0.3391
        \end{pmatrix},\\
&B = \begin{pmatrix}
      0.0843 & 0.0655\\
      0.4614 & 0.1750
    \end{pmatrix},
\bar{B} = \begin{pmatrix}
      0.8735 & 0.2435\\
      0.0124 & 0.7822
    \end{pmatrix},\\
&C = \begin{pmatrix}
      0.1499& 0.8148\\
0.1800 &0.7488
    \end{pmatrix},
\bar{C} = \begin{pmatrix}
      0.2084 &0.0877\\
0.6984 & 0.8266
    \end{pmatrix},\\
&D = \begin{pmatrix}
      0.9436 &0.0908\\
0.6800 & 0.1038
    \end{pmatrix},
\bar{D} = \begin{pmatrix}
      0.4319& 0.8377\\
0.3086 & 0.4565
    \end{pmatrix},\\
\end{aligned}
\end{equation*}
\textbf{weight coefficients:}
\begin{equation*}
\begin{aligned}
&Q = \begin{pmatrix}
     -0.0290 & 0.3912\\
0.3912 & 0.0257
    \end{pmatrix},
R = \begin{pmatrix}
      -0.1761 & 0.0362\\
0.0362 &0.5576
    \end{pmatrix},\\
&G = \begin{pmatrix}
      0.0940 & 0.2564\\
0.2564 &0.4806
    \end{pmatrix},\ \
\Gamma_1 = \begin{pmatrix}
      0.0815 &0.9821\\
0.1076 &0.6146
    \end{pmatrix}, \\
&\Gamma_2 = \begin{pmatrix}
      0.2807 &0.3933\\
0.0567 &0.1748
    \end{pmatrix}.
\end{aligned}
\end{equation*}
\begin{remark}
  Note that the eigenvalues of the weight matrices $Q$, $R$, $G$ are $\lambda_Q = (-0.3938, 0.3905)$, $\lambda_R = (-0.1779, 0.5594)$, $\lambda_G = (-0.0338, 0.6084)$ In this case, the standard assumption (i.e., $Q, G\geq 0$ and $R\gg0$), which has been widely applied in many previous research (see \cite{graber2016linear, yong2013linear, bensoussan2016linear, HCM2012}), fails to hold. However, through the regularity of Riccati equations (RE1)-(RE3), we can still design the control law for (\textbf{MC}) (\textbf{MG}) and (\textbf{MT}).
\end{remark}
Next, we verify the regularity of Riccati equations (RE1)-(RE3). The trajectories of the eigenvalues of $\mathcal{R}_1$, $\mathcal{R}_2$ and $\mathcal{R}_3$ are given as follows
\begin{figure}[H]
  \centering
  \includegraphics[width=\linewidth]{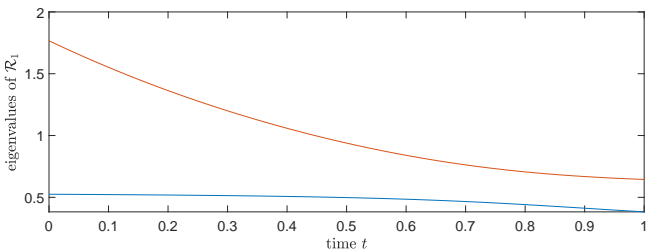}
  \caption{The eigenvalues of $\mathcal{R}_1$}\label{fig1}
\end{figure}
\begin{figure}[H]
  \centering
  \includegraphics[width=\linewidth]{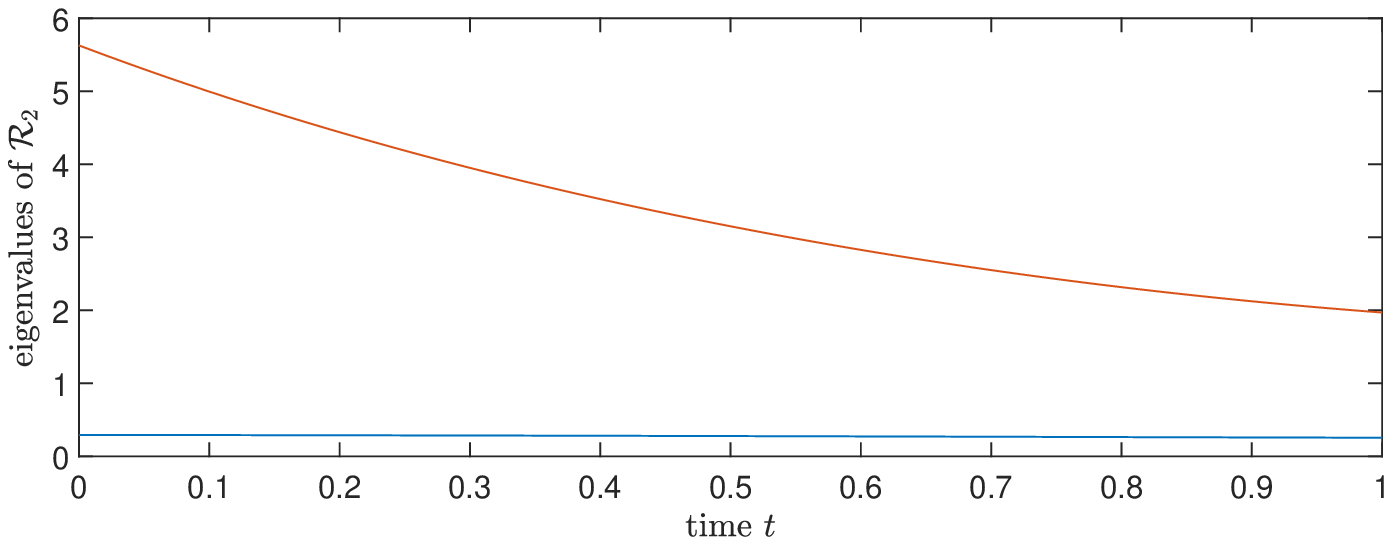}
  \caption{The eigenvalues of $\mathcal{R}_2$}\label{fig2}
\end{figure}
\begin{figure}[H]
  \centering
  \includegraphics[width=\linewidth]{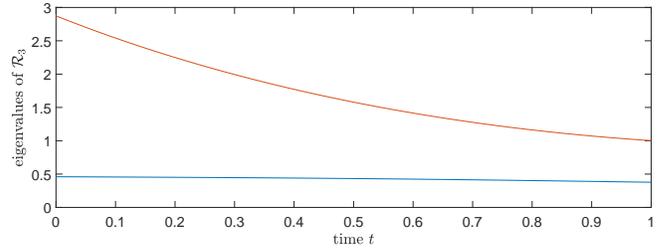}
  \caption{The eigenvalues of $\mathcal{R}_3$}\label{fig3}
\end{figure}
Through Fig. \ref{fig1}-\ref{fig3}, we see that $\mathcal{R}_1, \mathcal{R}_2, \mathcal{R}_3\gg 0$ on the whole time interval $[0,1]$. Next, we obtain the optimal control $\bar{u}_i$ and the mean field strategies $\tilde{\bu}^{\text{\rm \tiny MG}}$, $\tilde{\bu}^{\text{\rm \tiny MT}}$. Then the corresponding costs $\J_i^{\textrm{\tiny \rm MC}}(x_0;\bar{u}_i), \J_i^{\textrm{\tiny \rm MG}}(x_0;\tilde{\bu}^{\text{\rm \tiny MG}}), \J_i^{\textrm{\tiny \rm MT}}(x_0;\tilde{\bu}^{\text{\rm \tiny MT}})$ can be calculated. The trend of the corresponding costs with respect to the population $N$ is given by the following figure
\begin{figure}[H]
  \centering
  \includegraphics[width=\linewidth]{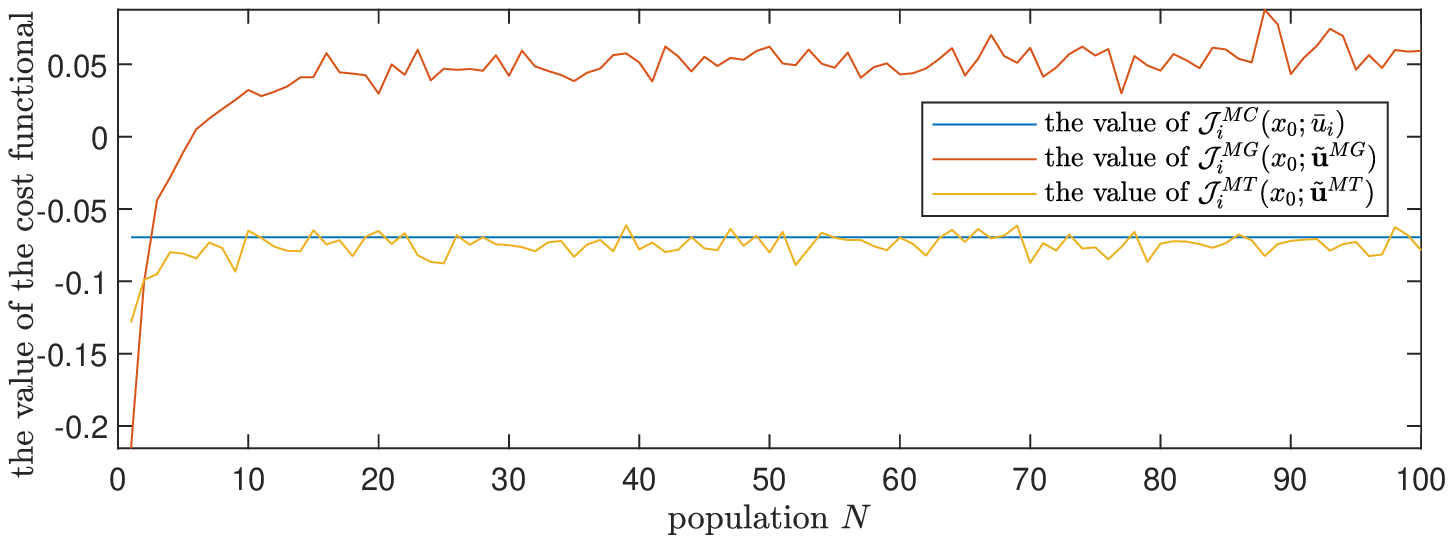}
  \caption{The trend of the costs w.r.t. the population}\label{fig4}
\end{figure}
Through Fig. \ref{fig4}, we see that $\J_i^{\textrm{\tiny \rm MT}}(x_0;\tilde{\bu}^{\text{\rm \tiny MT}} )\to \J_i^{\textrm{\tiny \rm MC}}(x_0;\bar{u}_i) $  and  $\J_i^{\textrm{\tiny \rm MG}}(x_0;\tilde{\bu}^{\text{\rm \tiny MG}} )> \J_i^{\textrm{\tiny \rm MC}}(x_0;\bar{u}_i)/ \J_i^{\textrm{\tiny \rm MT}}(x_0;\tilde{\bu}^{\text{\rm \tiny MT}} ) $ as $N\to \infty$, which is consistent with the result in Proposition \ref{prop4.8} (i.e., Fig. \ref{fig0}).

Moreover, by letting $\bar{A} = \bar{B} = \bar{C}  = \bar{D} = 0$, $\Gamma_1 \= \Gamma_2 \= I$ and the other coefficients  remaining unchanged, then we have the following trend of the corresponding costs with respect to the population $N$.
\begin{figure}[H]
  \centering
  \includegraphics[width=\linewidth]{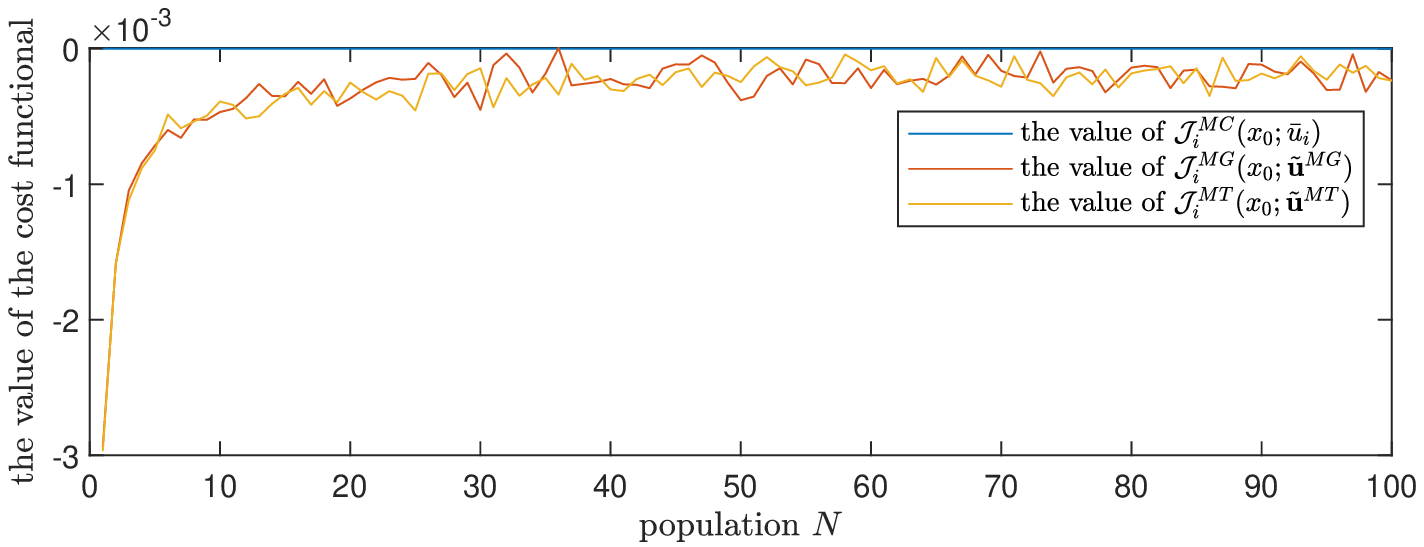}
  \caption{The trend of the costs w.r.t. the population}\label{fig5}
\end{figure}
Through figure \ref{fig5}, we see that
\[\J_i^{\textrm{\tiny \rm MT}}(x_0;\tilde{\bu}^{\text{\rm \tiny MT}} )\approx \J_i^{\textrm{\tiny \rm MG}}(x_0;\tilde{\bu}^{\text{\rm \tiny MG}} ) \to \J_i^{\textrm{\tiny \rm MC}}(x_0;\bar{u}_i), \]
 as $N\to \infty$, which is consistent with the result in Proposition \ref{prop4.8} (i.e., Fig. \ref{fig-1}).

\section{Conclusion} \label{sec 5}

This paper investigates a unified relation analysis of LQ MG, MT and MC. Based on stochastic maximum principle, fixed-point theory and person-by-person optimality principle, we derive the decentralized strategies, respectively. Two aspects of relations---decentralized strategies and value functions are obtained. Especially, two visualized figures (Fig. \ref{fig0} and Fig. \ref{fig-1}) are given to illustrate the relations of value functions, and a numerical example is provided to verify the theoretical results. It is worth pointing out that these theoretical results have abundant practical applications in engineering, economics, society science, etc. For example, when it is difficult or impossible to obtain the feasible strategy of complicated system (e.g. MT), we may study MG or MC as an \emph{alternative}, and the errors of value functions are also derived.

\end{document}